\documentclass[twoside,12pt]{article}
\usepackage{amsmath}
\usepackage{amssymb}
\usepackage[usenames,dvipsnames]{color}
\usepackage{cite}

\def\takeshi{\color{red}}
\def\pier{\color{blue}}

\def\colli{\color{Green}}
\def\fukao{\color{red}}
\let\takeshi\relax
\let\pier\relax

\let\colli\relax
\let\fukao\relax

\title{{C}ahn--{H}illiard equation with dynamic boundary conditions and mass constraint on the boundary}

\author{Pierluigi Colli\\
Dipartimento di Matematica, Universit\`a degli Studi di Pavia\\
Via Ferrata~1, 27100 Pavia, Italy\\
E-mail: \texttt{pierluigi.colli@unipv.it}\\
\and \\ Takeshi Fukao\\
Department of Mathematics, Faculty of Education\\
Kyoto University of Education\\
1~Fujinomori, Fukakusa, Fushimi-ku, Kyoto~612-8522 Japan\\
E-mail: \texttt{fukao@kyokyo-u.ac.jp}}

\newcommand\testopari{\sc Pierluigi Colli and Takeshi Fukao}
\newcommand\testodispari{\sc {C}ahn--{H}illiard with dynamic b.c.\ and mass constraint}
\begin{document}

\date{}

\maketitle

\begin{abstract}
The well-known {C}ahn--{H}illiard equation {\pier entails} mass conservation {\pier if
a suitable boundary condition is prescribed}. In the case when 
the equation {\pier {\takeshi is} also coupled with a} dynamic 
boundary condition, {\pier including} the {L}aplace--{B}e{\pier l}trami 
operator on the boundary, the total mass on the inside of the domain 
and its trace on the boundary should be conserved. 
The new issue of this paper is the setting 
of a mass constraint on the boundary. 
The effect of this additional constraint is 
{\pier the appearance of a} {L}agrange multiplier{\pier; in fact, 
two {L}agrange multipliers arise, one for 
the bulk, the other for} the boundary. {\pier 
The} well-posedness of the {\pier resulting} {C}ahn--{H}illiard {\pier system with} dynamic boundary condition 
and mass constraint on the boundary is obtained{\pier .} {\takeshi The theory of evolution} {\colli equations governed by {\pier subdifferentials} is exploited} and a complete
characterization of the solution {\takeshi is given}. 

\vspace{2mm}

\vspace{2mm}
\noindent \textbf{AMS (MOS) subject clas\-si\-fi\-ca\-tion:} 35K86, 49J40, 80A22.

\end{abstract}

%%%%% Section 1. %%%%%
\section{Introduction}

{\pier The famous {C}ahn--{H}illiard equation \cite{CH58, EZ86} offers a realistic description of the evolution 
phenomena related to solid-solid phase separation processes.
In this paper, we are interested to the mathematical investigation of it
and aim to analyze questions like existence and continuous dependence of solutions for  
a generalized {C}ahn--{H}illiard equation with dynamic boundary conditions 
and mass constraints on the boundary. Actually, we can solve the mathematical problem and, in particular, characterize the constraint with the help of a {L}agrange multiplier.}

Let $0<T<+\infty$ and {\pier let} $\Omega \subset \mathbb{R}^{d}$, $d=2$ or $3$, be the bounded smooth 
domain occupied by the material.  Also the boundary $\Gamma$ of $\Omega $ is supposed to be smooth enough. 
We recall the isothermal {C}ahn--{H}illiard equation in the following generalized form:
\begin{gather*} 
	\frac{\partial u}{\partial t}-\Delta \mu = 0 
	\quad \mbox{in }Q:=\Omega \times (0,T), 
	\\
	 {\pier \mu = \tau \frac{\partial u}{\partial t}-\Delta u + \xi +\pi (u)-f ,
	\quad \xi \in \beta (u)} \quad\mbox{in }Q,
\end{gather*} 
where {\pier the} unknowns $u:=u(x,t)$ {\pier and} $\mu :=\mu (x,t)$ stand for the order parameter {\pier and} the chemical potential, respectively. 
Moreover, {\pier $\tau $ is a viscosity coefficient which can be greater or equal to $0$ (we treat both cases);}
$\beta $ stands for the {\pier subdifferential} of the convex part $\widehat{\beta }$ and 
$\pi $ stands for {\pier the derivative of} 
the concave perturbation $\widehat{\pi}$ of a 
double well potential {\pier $W= \widehat{\beta }+ \widehat{\pi }$}, 
for example {\pier $W (r) = (r^2 -1)^2/4 $ with $\beta (r)=r^3$ and 
$\pi (r)=-r$} for all  $r \in \mathbb{R}$. 
{\pier In general, $\beta $ is assumed to be a} 
maximal monotone graph in $\mathbb{R} \times \mathbb{R}$. 
Recently, {\pier this equation was treated in some papers \cite{CGS14, CGS14bis, GMS09, GMS10} 
when coupled with a dynamic boundary condition of the following form:}
\begin{gather*} 
	u_\Gamma = u_{|_\Gamma }
	\quad \mbox{on }\Sigma:={\pier \Gamma}\times (0,T),
	\\
	{\pier \partial _\nu u
	+ \frac{\partial u_\Gamma }{\partial t}-\Delta _{\Gamma } u_\Gamma 
	+ \xi _\Gamma +\pi _\Gamma (u_\Gamma ) = f_\Gamma ,
	\quad 	\xi _\Gamma \in \beta _\Gamma (u_\Gamma) }
	\quad \mbox{on }\Sigma,
\end{gather*} 
where, $u_{|_\Gamma }$ denotes the trace of $u$ and 
$\partial _\nu $ represents the outward normal derivative on $\Gamma $. 
$\Delta _{\Gamma }$ stands for the {L}aplace--{B}e{\pier l}trami operator
on $\Gamma $ (see, e.g., \cite[Chapter 3]{Gri09}), 
$\beta _\Gamma$ and $\pi _\Gamma $ have the same property as $\beta $ and $\pi$, respectively. 

About dynamic boundary conditions, {\pier let us point out that}
the mathematical research for the various problem was already running 
in {\pier the} 1990's. For example, the {S}tefan problem with dynamic boundary conditions  
was treated in the series of {A}iki \cite{Aik93, Aik95, Aik96}. 
Recent advances in the {C}ahn--{H}illiard equation with 
the dynamic boundary conditions
can be found in 
\cite{CGS14, GMS09, GMS10, GM13, RZ03} and references therein. 

{\pier As is well known}, conservation of $u$ is required. 
Therefore, under the homogeneous {N}eumann boundary 
condition 
$$
	\partial _\nu \mu  = 0 
	\quad \mbox{on }\Sigma,
$$
we can realize that 
$$
	\frac{1}{|\Omega |} \int_{\Omega }^{} u(t) dx =m_0:=\frac{1}{|\Omega |} \int_{\Omega }^{} u_0 dx 
	\quad \mbox{for all } t \in [0,T],
$$
for {\pier a} given initial data $u_0$. 
The new issue of this paper is the setting of {\pier a}  mass constraint on the 
boundary. 
More precisely, we require that the solution $u$ satisfies
$$
	k_* \le 
	\int_{\Gamma }^{} w_{\Gamma }{\pier u_\Gamma} (t) d\Gamma 
	\le k^* 
	\quad \mbox{for all }t \in [0,T],
$$
where {\pier $k_*$ and $k^*$ are fixed constants fulfilling $k_* \le k^*$ and}
$w_{\Gamma }$ is given weight function on $\Gamma$. 
This kind of problem for the {A}llen--{C}ahn equation 
was treated {\pier in} \cite{CF14}, by applying 
{\pier the abstract theory developed in} \cite{FK13}.  
{\pier In the case of the {C}ahn--{H}illiard equation, the 
essential structure of the constraint has been studied in} 
\cite{KN96, Kub12}. 
We can also find {\pier a similar treatment for 
the preservation of the constraint} in \cite{Aik96, CGM13}. 

A brief outline of the present paper {\pier along with 
a short description of the various items is as follows.}

In Section 2, we present the main results, 
consisting in the well-posedness of the {C}ahn--{H}illiard equation with 
dynamic boundary conditions and mass constraints on the boundary. 
We write the system as an 
evolution inclusion and characterize the solution with the help of the 
{L}agrange multipliers.
We also remark {\pier that actually there will be} two {L}agrange multipliers. 

In Section 3, we prove the continuous dependence {\pier and
of course this result} entails the uniqueness property. 

In Section 4, we prove the existence result. 
{\pier The proof is split in several steps. 
First, we construct an approximate solution 
by substituting the maximal monotone graphs with their 
{M}oreau--{Y}osida regularizations, in the case when $\tau >0$. 
The solvability of the approximate problem is guaranteed by the 
abstract theory of doubly nonlinear evolution inclusions~\cite{CV90}. 
Moreover, arguing in a similar way as in \cite{FK13},
we show that the solution satisfies suitable regularity properties  
and obtain a strong characterization of the approximate problem 
by the Lagrange multiplier: in fact, we are able to prove uniform a priori estimates
on all the components of the solution. 
And finally, from these estimates, 
we can pass to the limit and conclude the existence proof in the case $\tau >0$. 
Next, we can proceed by considering the limiting problem as $\tau \to 0$ and
derive the well-posedness result in the pure {C}ahn--{H}illiard case as well.}

\begin{itemize}
 \item[1.] Introduction
 \item[2.] Main results
\begin{itemize}
 \item[2.1.] Definition of the solution by the {L}agrange multiplier
 \item[2.2.] Remark for the {L}agrange multiplier
 \item[2.3.] Well-posedness 
 \item[2.4.] Abstract formulation 
\end{itemize}
 \item[3.] Continuous dependence 
 \item[4.] Existence
\begin{itemize}
 \item[4.1.] Approximation of the problem
 \item[4.2.] A priori estimates
 \item[4.3.] Passage to the limit as $\varepsilon \to 0$
 \item[4.4.] Passage to the limit as $\tau \to 0$
\end{itemize}
\end{itemize}

%%%%% Section 2. %%%%%
\section{Main results}

In this section, we present our main result,  
{\pier which states} the well-posedness of the {C}ahn--{H}illiard 
equation with the dynamic boundary conditions 
and mass constraints on the boundary. 
We apply the treatment of the dynamic boundary conditions 
{\pier as in} \cite{CC13, CF14} and {\pier exploit} the abstract 
theory of the evolution inclusion, essentially {\pier referring} 
to \cite{FK13, KN96}.

%%%%% Section 2.1. %%%%%
\subsection{Definition of the solution by the {L}agrange multiplier}

Let $0<T<+\infty$ and $\Omega \subset \mathbb{R}^d$, 
$d=2$ or $3$, be the bounded domain 
with smooth boundary $\Gamma:=\partial \Omega $. 
We use the notation:
\begin{align*}
&H_0:=L^2(\Omega)_0:={\pier \Bigl\{} z \in L^2(\Omega )\ : \  \int_{\Omega }^{} z dx =0 {\pier \Bigr\}}, \\[0.1cm]
&{\pier H_\Gamma :=L^2(\Gamma), \quad
V_0:=H^1(\Omega) \cap H_0, \quad V_{\Gamma }:=H^1(\Gamma ),}
\end{align*}
with usual norms {\pier%
$|\cdot |_{H_0}$, 
$|\cdot |_{H_{\Gamma}}$, 
$$
|z|_{V_0}:=|\nabla z|_{L^2(\Omega )^d} \ \hbox{ for } \ z \in V_0, \ \quad  
|z_\Gamma |_{V_{\Gamma }} := {\takeshi 
\left\{  \int_\Gamma \left( |z_\Gamma|^2 + |\nabla_\Gamma 
z_\Gamma|^2 \right) d\Gamma \right\}^{\frac{1}{2}}} \ \hbox{ for } \ z_\Gamma \in V_\Gamma,
$$ 
respectively. Here, $\nabla _{\Gamma }$ denotes the surface gradient on $\Gamma$ (see, e.g., \cite[Chapter 3]{Gri09}). Moreover, let $V_0^*$ be} the dual space of $V_0$ and 
$F:V_0 \to V_0^*$ {\pier denote} the duality mapping defined by 
$$
	\langle Fy , z \rangle _{V^*_0,V_0} := \int_{\Omega }^{} \nabla y \cdot \nabla z dx 
	\quad \mbox{for all } \, y,z \in V_0.
$$
Then, {\pier the form $(\cdot ,\cdot )_{V_0^*}: V_0^*\times V_0^* \to \mathbb{R}$,}
$$
	(y^* ,z^* )_{V_0^*}:=\int_{\Omega }^{} \nabla F^{-1}y^* \cdot \nabla F^{-1}z^* dx 
	\quad \mbox{for all } \, y^*,z^* \in V_0^*,
$$
{\pier yields} the inner product {\pier in}  $V_0^*$. {\pier Here, 
$F^{-1}$ is the inverse operator of $F$ and its restriction to $H_0$ works as follows: 
if $z \in H_0$, $y=F^{-1}z$ uniquely solves the boundary value problem
$$
	\left\{ 
	\begin{array}{ccl}
	-\Delta y & = & z \quad \mbox{a.e. in }\Omega, \\
	\displaystyle \partial_\nu y & = & 0 \quad \mbox{a.e. on } \Gamma, \vspace{2mm}\\
	\displaystyle \int_{\Omega }^{} y dx & = & 0.
	\end{array} 
	\right. 
$$
and consequently lies in $H^2(\Omega )$, due to well-known elliptic regularity results.} The reader can check that testing $-\Delta y = z$ by some $\tilde{z} \in V_0$
{\pier leads to}
$$ 
	\int_{\Omega }^{} \nabla y \cdot \nabla \tilde{z} dx = \int_{\Omega }^{} z \tilde{z} dx \quad 
	\mbox{for all } \tilde{z} \in V_0,
$$
that is, $z=Fy$ {\pier as expected}. {\pier Finally, by virtue of the 
{P}oincar\'e--{W}irtinger inequality} there exists 
a constant ${\pier C_0}>0$ such that 
\begin{equation}
{\pier |z|_{H_0}^2\le {\pier C_0} |z|_{V_0}^2 \ \hbox{ for all }\, z \in V_0. } \label{PW}
\end{equation} 
Then, we obtain
$V_0 \mathop{\hookrightarrow} \mathop{\hookrightarrow} 
H_0 \mathop{\hookrightarrow} \mathop{\hookrightarrow}V_0^*$, where 
``$\mathop{\hookrightarrow} \mathop{\hookrightarrow} $'' stands for 
the dense and compact embedding, namely 
$(V_0,H_0,V_0^*)$ is {\pier a} standard {H}ilbert triplet. 
The same considerations hold for $H_{\Gamma }$ and 
$V_{\Gamma }$. {\pier Now, we set}
$$
	\mbox{\boldmath $ H $}_0:=
	H_0
	\times
	H_{\Gamma },
	\quad 
	\mbox{\boldmath $ V $}_0:=
	\bigl\{ 
	(u,u_{\Gamma }) \in 
	V_0 \times V_{\Gamma }\,: \ 
	u_{|_\Gamma }=u_{\Gamma }~ \, \mbox{a.e.\ on } \Gamma 
	\bigr\},
$$
where $u_{|_\Gamma}$ denotes the trace of $u$. Observe that 
$\mbox{\boldmath $ H $}_0$ and $\mbox{\boldmath $ V $}_0$ are {H}ilbert spaces 
with the inner products 
\begin{align*}
	(\mbox{\boldmath $ u $},\mbox{\boldmath $ z $})_{\mbox{\scriptsize \boldmath $ H$}_0}
	:=(u,z)_{H_0} + (u_\Gamma ,z_\Gamma )_{H_\Gamma } \quad 
	\mbox{for all}~\mbox{\boldmath $ u$}
	:=(u,u_{\Gamma }), \,
	\mbox{\boldmath $ z$}:=(z,z_{\Gamma }) 
	\in \mbox{\boldmath $ H$}_0,\\
	(\mbox{\boldmath $ u $},\mbox{\boldmath $ z $})_{\mbox{\scriptsize \boldmath $ V$}_0}
	:=(u,z)_{V_0} + (u_\Gamma ,z_\Gamma )_{V_\Gamma } \quad 
	\mbox{for all}~\mbox{\boldmath $ u$}
	:=(u,u_{\Gamma }), \,
	\mbox{\boldmath $ z$}:=(z,z_{\Gamma }) 
	\in \mbox{\boldmath $ V$}_0
\end{align*}
{\pier and related norms.}
Then, we obtain 
$\mbox{\boldmath $ V$}_0 
\mathop{\hookrightarrow} \mathop{\hookrightarrow}
\mbox{\boldmath $ H$}_0 
\mathop{\hookrightarrow} \mathop{\hookrightarrow}
\mbox{\boldmath $ V$}_0^*$ (see, e.g., \cite[Appendix]{CF14}).
As a remark, let us restate that if 
$\mbox{\boldmath $ u$}=(u,u_{\Gamma}) \in \mbox{\boldmath $ V$}_0$ then $u_{\Gamma }$ 
is exactly the trace of $u$ on $\Gamma$, while, if $\mbox{\boldmath $ u$}=(u,u_{\Gamma})$ is just in $ \mbox{\boldmath $ H$}_0$, then $u \in H$ 
and $u_{\Gamma } \in H_{\Gamma }$ are independent.

The initial-value problem for the {C}ahn--{H}illiard equation with 
dynamic boundary conditions {\pier can be set} as the following system \eqref{(1)}--\eqref{(6)}: 
\begin{gather} 
	\frac{\partial u}{\partial t}-\Delta \mu =0
	\quad \mbox{in } Q,
	\label{(1)}
	\\ 
	{\pier \mu = \tau \frac{\partial u}{\partial t}-\Delta u 
	+ \xi +\pi (u)-f ,
	\quad \xi \in \beta (u)}
	\quad \mbox{in }Q,
	\label{(2)}
	\\ 
	\partial _\nu \mu =0
	\quad \mbox{on } \Sigma,
	\label{(3)}
	\\
	u_\Gamma=u_{|_\Gamma }, 
	\quad \mbox{on } \Sigma,
	\label{(4)}
	\\ 
	{\pier 
	\partial _\nu u 
	+ \frac{\partial u_\Gamma }{\partial t}-\Delta _{\Gamma } u_\Gamma
	+ \xi _{\Gamma }+\pi _{\Gamma }(u_\Gamma) = f_\Gamma ,\quad \xi _\Gamma \in \beta _\Gamma (u_\Gamma)}
	\quad \mbox{on } \Sigma,
	\label{(5)}
	\\
	u(0)=u_0 \quad \mbox{in } \Omega, 
	\quad u_\Gamma (0)=u_{0\Gamma} \quad \mbox{on } \Gamma,
	\label{(6)}
\end{gather}
where $\tau \ge 0$ is a viscosity coefficient. {\pier Testing \eqref{(1)} by the constant function $1$ and using the boundary condition \eqref{(3)}, we realize that 
{\takeshi $\partial u/\partial t $} has zero mean value {\takeshi in $\Omega$.} Then, a formal test of \eqref{(1)} and \eqref{(2)} by an arbitrary element $z\in V_0$ and a subsequent combination produce, with the help of the definition 
of $F$ and the conditions in \eqref{(3)}--\eqref{(5)}, the variational formulation
\begin{align}
	\lefteqn{ 
	\int_{\Omega }^{} F^{-1} \left( \frac{\partial u}{\partial t}(t)\right)z dx 
	+ \tau \int_{\Omega }^{} \frac{\partial u}{\partial t}(t) z dx 
	+ \int_{\Gamma }^{} \frac{\partial u_\Gamma }{\partial t}(t) z_\Gamma d \Gamma 
	} \nonumber \\
	&  \quad {} 
	+ \int_{\Omega }^{} \nabla u(t) \cdot \nabla z dx 
	+ \int_{\Gamma }^{} \nabla_\Gamma u_\Gamma(t) \cdot \nabla_\Gamma z_\Gamma d\Gamma 
	+ \int_{\Omega }^{} \xi(t) z dx 
	+ \int_{\Gamma }^{} \xi _\Gamma(t) z_\Gamma d\Gamma 
	\nonumber \\
	& \quad 
	{}
	+ \int_{\Omega }^{} \pi \bigl( u(t) \bigr) z dx 
	+ \int_{\Gamma }^{} \pi _\Gamma \bigl(u_\Gamma (t)\bigr) 
	z_\Gamma d\Gamma 
	 =  \int_{\Omega }^{} f(t) z dx 
	+ \int_{\Gamma }^{} f_\Gamma(t) z_\Gamma d\Gamma,
	\label{pier1}
\end{align}
for a.e.\ $t\in (0,T)$, 
for all $z \in V_0$ with {\takeshi $z_\Gamma = z_{|_\Gamma} $}. We are now interested to deal not directly with \eqref{pier1} but with a variational inequality replacing it, where the solution and the test function vary in a suitable convex set.}

{\pier Concerning the data,} we assume that
\begin{itemize}
 \item[(A1)] {\pier $\beta $, $\beta _{\Gamma }$, maximal monotone graphs in 
$\mathbb{R} \times \mathbb{R}$, are the subdifferentials} 
$$
	\beta =\partial \widehat{\beta}, \quad \beta _{\Gamma }=\partial \widehat{\beta }_{\Gamma }
$$
of some {\pier continuous} and convex functions 
$${\pier \widehat{\beta }, \, \widehat{\beta }_{\Gamma }: \mathbb{R} \to [0,+\infty ) \ \hbox{ such that } \ 
\widehat{\beta }(0)=\widehat{\beta}_{\Gamma }(0)=0;}$$ 
 \item[(A2)]
{\pier  
$\pi $, $\pi _{\Gamma }: \mathbb{R} \to \mathbb{R}$ are {L}ipschitz continuous functions}
with {L}ipschitz constants $L$ and $L_{\Gamma}$, respectively;
 \item[(A3)] $\mbox{\boldmath $ f$} := (f,f_\Gamma ) \in L^2(0,T;L^2(\Omega )) \times L^2(0,T;H_\Gamma )$ and 
$\mbox{\boldmath $ u$}_0 := (u_0,u_{0\Gamma }) \in H^1(\Omega ) \times V_\Gamma ${\pier , where $u_{0\Gamma }:= {u_0}_{|_\Gamma }$.}
\end{itemize} 
In particular, {\pier by (A1) we are asking that} 
$D(\beta )=D(\beta _{\Gamma })=\mathbb{R}$, 
$0 \in \beta (0)$ and $0 \in \beta _{\Gamma }(0)$.

In this paper, we are interested to the setting of {\pier the constraint}
\begin{equation} 
	k_* \le \int_{\Gamma }^{} w_{\Gamma }u_{|_\Gamma }(t) d\Gamma 
	\le k^* 
	\quad \mbox{for all }t \in [0,T],
	\label{const.}
\end{equation} 
{\pier for the solution to the related variational inequality 
(cf.\ \eqref{pier1}).}
Here, $k_*$ and $k^*$ are real constants with $k_* \le k^*$, and 
$\mbox{\boldmath $ w$}:=(0,w_{\Gamma }) \in \mbox{\boldmath $ H$}_0$ 
is fixed. We require that the weight function $w_\Gamma $ 
satisfies 
\begin{itemize}
 \item[(A4)] {\takeshi $w_\Gamma\in H_\Gamma$, $w_{\Gamma } \ge 0$ a.e.\ on 
$ \Gamma$ and $\sigma_0 :=\int_{\Gamma }^{} w_{\Gamma } d\Gamma >0$. }
\end{itemize}
The last inequality can be seen as a nondegeneracy
condition on the weight element $\mbox{\boldmath $ w$}$. 
\vspace{2mm}

{\takeshi Hence}, let us term {\rm (P)} the initial-value problem related to the variational inequality
and to the constraint in \eqref{const.}. Now, we define precisely the notion of solution to the
problem {\rm (P)} by means of a {L}agrange multiplier. 
In order to {\pier set $H_0$ as the pivot} space, put 
$m_0:=(1/|\Omega |) \int_{\Omega }^{} u_0 dx$ and 
let $v(x,t):=u(x,t)-m_0$ be the new unknown function and {\pier define} analogously
$v_0:=u_0-m_0$ in $\pier \Omega$, $v_{0\Gamma }:=u_{0\Gamma }-m_0$ on 
$\pier\Gamma$, $h_*:=k_*-m_0 \sigma _0$ and 
$h^*:=k^*-m_0 \sigma _0$, respectively.

\paragraph{Definition 2.1.} 
{\it The quadruplet $(\mbox{\boldmath $ v$}, \mbox{\boldmath $ {\pier\xi} $}, \omega, \lambda )$ 
is called the solution of {\rm (P)} if
\begin{align*}
	\mbox{\boldmath $ v$}=(v,v_{\Gamma }) \quad \mbox{with} \quad
	v \in H^1(0,T;H_0) \cap C \bigl( [0,T];V_0 \bigr) \cap L^2 \bigl( 0,T;H^2(\Omega ) \bigr),\\
	v_{\Gamma } \in H^1(0,T;H_{\Gamma}) \cap C \bigl( [0,T];V_{\Gamma } \bigr) \cap L^2 \bigl( 0,T;H^2(\Gamma ) \bigr),\\
	\mbox{{\boldmath ${\pier\xi}$}}=({\pier\xi} ,{\pier\xi}_{\Gamma }) \in  L^2(0,T;\mbox{\boldmath $H$}_0) , \quad
	\omega, \lambda \in L^2(0,T),
\end{align*}
and $v$, $v_\Gamma$, ${\pier\xi}$, ${\pier\xi}_\Gamma$, $\omega $, $\lambda $ satisfy
\begin{gather} 
	F^{-1}\left( \frac{\partial v}{\partial t} \right) 
	+ \tau \frac{\partial v}{\partial t}-\Delta v + {\pier\xi} + \pi(v+m_0)= f + \omega 
	\quad \mbox{a.e.\ in } Q,
	\label{(D1)}\\
	{\pier\xi} \in \beta (v+m_0)
	\quad \mbox{a.e.\ in } Q,
	\label{(D2)}\\
	v_\Gamma=v_{|_\Gamma }, \quad 
	\partial_\nu v 
	+ \frac{\partial v_{\Gamma }}{\partial t}-\Delta _{\Gamma } v_{\Gamma }
	+ {\pier\xi} _{\Gamma } +\pi _{\Gamma }(v_{\Gamma }+m_0) + \lambda w _{\Gamma } = f_\Gamma 
	\quad \mbox{a.e.\ on } \Sigma,
	\label{(D3)}\\ 
	{\pier\xi}_\Gamma \in \beta_{\Gamma } (v_{\Gamma }+m_0)
	\quad \mbox{a.e.\ on } \Sigma,
	\label{(D4)}\\
	v(0)=v_0 \quad \mbox{a.e.\ in } \Omega, \quad 
	v_{\Gamma }(0)=v_{0\Gamma} \quad \mbox{a.e.\ on } \Gamma,
	\label{(D5)}\\
	h_* \le \int_{\Gamma }^{} w_{\Gamma } v_{\Gamma }(t)d\Gamma 
	\le h^*
	\quad \mbox{for a.a.\ } t \in (0,T),
	\label{(D6)}
	\\ {\pier \lambda (t) \int_{\Gamma }^{} w_{\Gamma } \bigl( v_{\Gamma }(t) - z_\Gamma 
 \bigr) d\Gamma  \geq 0 \nonumber \quad  \hbox{for a.a. } t \in (0,T)\qquad \quad \ \qquad}\\ 
{\pier \mbox{and for all\ }  \mbox{\boldmath $z$} =(z, z_{\Gamma }) 
\in \mbox{\boldmath $V$}_0
\  \mbox{such that } \ h_* \le  \int_{\Gamma }^{} w_{\Gamma } z_{\Gamma }\, d\Gamma \le h^*.}
 \label{pier5}
\end{gather} 
In the case $\tau =0$, the regularity of $v$ should be modified into} 
$$
	v \in H^1(0,T;V_0^*) \cap L^\infty (0,T;V_0) \cap L^2 \bigl (0,T;H^2(\Omega ) \bigr ).
$$

%%%%% Section 2.2. %%%%%
\subsection{Remark for the {L}agrange {\pier multipliers}}

By comparing \eqref{(2)} with \eqref{(D1)}--\eqref{(D2)}, we realize that 
$$
	\mu =- F^{-1}\left( \frac{\partial v}{\partial t}\right) + \omega \quad \mbox{a.e.\ in } Q,
$$
{\pier so that $\omega $ turns out to be} the mean value of the chemical potential $\mu $
$$
	\omega (t)=\frac{1}{|\Omega |} \int_{\Omega }^{} \mu (t) dx.
$$
{\pier On the other hand, $\lambda $} has the role of a {L}agrange multiplier related to the constraint in \eqref{(D6)} 
on the boundary. Then, the two {L}agrange multipliers $\omega $ and 
$\lambda $ have different meaning{\pier ; in particular,} $\lambda $ is obtained by 
solving the problem and it explicitly appears in the variational formulation, while $\omega $ does not {\pier show up in the variational 
inequality and it can be only} identified a posteriori. 
Indeed, if we test \eqref{(D1)} by a function $z \in V_0$, 
then $\omega $ disappears %from the variational equality 
and we obtain
{\pier (cf.\ also \eqref{pier1})}
\begin{align}
	\lefteqn{ 
	\int_{\Omega }^{} F^{-1} \left( \frac{\partial v}{\partial t}(t)\right)z dx 
	+ \tau \int_{\Omega }^{} \frac{\partial v}{\partial t}(t) z dx 
	+ \int_{\Gamma }^{} \frac{\partial v_\Gamma }{\partial t}(t) z_\Gamma d \Gamma 
	} \nonumber \\
	&  \quad {} 
	+ \int_{\Omega }^{} \nabla v(t) \cdot \nabla z dx 
	+ \int_{\Gamma }^{} \nabla_\Gamma v_\Gamma(t) \cdot \nabla_\Gamma z_\Gamma d\Gamma 
	+ \int_{\Omega }^{} {\pier\xi}(t) z dx 
	+ \int_{\Gamma }^{} {\pier\xi} _\Gamma(t) z_\Gamma d\Gamma 
	\nonumber \\
	& \quad 
	{}
	+ \int_{\Omega }^{} \pi \bigl( v(t)+m_0 \bigr) z dx 
	+ \int_{\Gamma }^{} \pi _\Gamma \bigl(v_\Gamma (t)+m_0 \bigr) z_\Gamma d\Gamma 
	+ \int_{\Gamma }^{} \lambda (t)w_\Gamma z_\Gamma d\Gamma 
	\nonumber \\
	& =  \int_{\Omega }^{} f(t) z dx 
	+ \int_{\Gamma }^{} f_\Gamma(t) z_\Gamma d\Gamma,
	\label{v.i.}
\end{align}
for all $z \in V_0$ satisfying $z_{|_\Gamma} =z_\Gamma$,
because $(\omega (t),z)_{H_0}=0$. 
On the contrary, if we simply integrate {\pier \eqref{(D1)} and set}
\begin{equation}
\label{pier2}
	q := {\pier\xi} +\pi (v+m_0)-f 
	\quad {\rm a.e.\ in}~Q, \quad 
	q_\Gamma :={\pier\xi} _\Gamma +\pi_\Gamma (v_\Gamma +m_0)-f_\Gamma
	\quad {\rm a.e.\ on}~\Sigma,
\end{equation}
{\pier with the help of \eqref{(D3)}} we obtain 
\begin{equation} 
	\omega (t)= \frac{1}{|\Omega |} 
	\left\{ \int_{\Omega }^{}  q(t) dx + 
	\int_{\Gamma }^{} \left( \frac{\partial v_\Gamma }{\partial t}(t) + q_\Gamma (t)+\lambda (t)w_\Gamma 
	\right) d\Gamma \right\}
	\quad \mbox{for all } t \in [0,T]. \label{omega}
\end{equation} 

In the last part of this section, we {\pier show how to} 
recover \eqref{(D1)} and \eqref{(D3)} from the variational equality \eqref{v.i.}. 
Define the projection $P_0 :L^2(\Omega ) \to H_0$ by 
$$
	P_0 z :=z - \frac{1}{|\Omega |} \int_{\Omega }^{} z dx \quad {\rm for~all}~z \in L^2(\Omega ). 
$$
{\pier Take $z \in H_0^1(\Omega )$ (so that $z_{|_\Gamma} =0 $ {\takeshi a.e.\ on} $\Gamma$) and use $P_0z$ as test function in \eqref{v.i.}. We note that 
$(P_0z)_{|_\Gamma} =-(1/|\Omega |)\int_{\Omega }^{} z dx$ and infer}\begin{align*}
	\lefteqn{ 
	\int_{\Omega }^{} F^{-1} \left( \frac{\partial v}{\partial t}(t)\right)z dx 
	+ \tau \int_{\Omega }^{} \frac{\partial v}{\partial t}(t) z dx 
	+ \int_{\Gamma }^{} \frac{\partial v_\Gamma }{\partial t}(t) d \Gamma \left( - \frac{1}{|\Omega |} \int_{\Omega }^{}z d\tilde{x} \right)
	} \nonumber \\
	& \quad {} 
	+ \int_{\Omega }^{} \nabla v(t) \cdot \nabla z dx 
	+ \int_{\Omega }^{} \Bigl( {\pier\xi}(t)+ \pi \bigl( v(t)+m_0 \bigr)-f(t) \Bigr) \left( z - \frac{1}{|\Omega |} \int_{\Omega }^{}z d\tilde{x} \right)  dx 
	\nonumber \\
	& \quad {}
	+ \int_{\Gamma }^{} \Bigl( {\pier\xi} _\Gamma(t) + \pi _\Gamma \bigl(v_\Gamma (t)+m_0\bigr) -f_\Gamma(t) \Bigr) d\Gamma 
	\left( - \frac{1}{|\Omega |} \int_{\Omega }^{}z d\tilde{x} \right) 
	{\pier {}= 0.} 
\end{align*}
Then, {\pier recalling the notation \eqref{pier2} we {\pier easily} obtain the equation in the interior, i.e.,} 
$$
	F^{-1} \left( \frac{\partial v}{\partial t}\right)
	+ \tau \frac{\partial v}{\partial t}
	- \Delta v
	+ P_0 q
	- \frac{1}{|\Omega |} 
	\int_{\Gamma }^{} 
	\left( \frac{\partial v_\Gamma }{\partial t} + q_\Gamma  
	\right) d\Gamma
	= 0 \quad \mbox{a.e.\ in } Q
$$
{\pier and, in view of \eqref{omega}, we find out that} 
$$
	F^{-1}\left( \frac{\partial v}{\partial t} \right) 
	+ \tau \frac{\partial v}{\partial t}-\Delta v + {\pier q} =  \omega
	\quad \mbox{a.e.\ in } Q.
$$
Next, we take a general $\mbox{\boldmath $ z$}:=(z,z_\Gamma) \in \mbox{\boldmath $ V$}_0$ and 
note that {\pier \eqref{v.i.} reduces to} 
\begin{align*}
	%\lefteqn{ 
	&\int_{\Omega }^{} F^{-1} \left( \frac{\partial v}{\partial t}(t)\right)\! z dx 
	+ \tau \int_{\Omega }^{} \frac{\partial v}{\partial t}(t) z dx 
	+ \int_{\Gamma }^{} \frac{\partial v_\Gamma }{\partial t}(t) z_\Gamma d \Gamma 
	- \int_{\Omega }^{} \Delta v(t) z dx 
	{\pier {} + \int_{\Gamma }^{} \partial_\nu v(t) z_\Gamma d \Gamma }
	% } 
	\nonumber \\
	& \quad {} 
	+ \int_{\Gamma }^{} \nabla_\Gamma v_\Gamma(t) \cdot \nabla_\Gamma z_\Gamma d\Gamma 
	+ \int_{\Omega }^{} q(t) z dx 
	+ \int_{\Gamma }^{} q_\Gamma(t) z_\Gamma d\Gamma 
	+ \int_{\Gamma }^{} \lambda (t)w_\Gamma z_\Gamma d\Gamma 
%	\nonumber \\ & 
	{\pier {}= 0,} 
\end{align*}
which means that 
$$
	\int_{\Omega }^{} \omega(t) z dx + \int_{\Gamma }^{} 
	\left( \partial _\nu v(t) 
	+\frac{\partial v_\Gamma }{\partial t}(t)
	-\Delta_\Gamma v_\Gamma (t)+q_\Gamma (t)+\lambda (t) w_\Gamma  
	\right) z_\Gamma d\Gamma =0.
$$
By virtue of the fact {\pier that $\int_{\Omega }^{} \omega(t) z dx=\omega (t)\int_{\Omega }^{}  z dx=0$, we finally have (cf.\ \eqref{(D3)})}
$$
	\partial _\nu v
	+ \frac{\partial v_{\Gamma }}{\partial t}-\Delta _{\Gamma } v_{\Gamma }
	+ {\pier q _{\Gamma } + \lambda w _{\Gamma } = 0}
	\quad \mbox{a.e.\ on } \Sigma.
$$

%%%%% Section 2.3. %%%%%
\subsection{Well-posedness}

The first result states the continuous dependence on the data. 
The uniqueness of the component $\mbox{\boldmath $ v$}$ of 
the solution is also guaranteed by this theorem. 

\paragraph{Theorem 2.1.} 
{\it Let $\tau \ge 0$. Assume {\rm (A1)}--{\rm {\pier (A4)}}. For $i=1,2$, 
let $( \mbox{\boldmath $ v$}^{(i)}, \mbox{\boldmath $ {\pier\xi} $}^{(i)}, \omega ^{(i)}, \lambda^{(i)} )$, 
with $\mbox{\boldmath $ v$}^{(i)}=(v^{(i)},v^{(i)}_\Gamma )$ and 
$\mbox{\boldmath $ {\pier\xi} $}^{(i)}=({\pier\xi} ^{(i)},{\pier\xi} ^{(i)}_\Gamma)$
be a solution to {\rm (P)} corresponding to the data 
$\mbox{\boldmath $ f$}^{(i)}=(f^{(i)}, f^{(i)}_{\Gamma })$ and 
$\mbox{\boldmath $ v$}_0^{(i)}=(v_0^{(i)}, v_{0\Gamma }^{(i)})$. 
Then, there 
exists a positive constant $C>0$, 
depending on $L$, $L_{\Gamma}$ and $T$, such that} 
\begin{align} 
	\lefteqn{ 
	\bigl| v^{(1)}(t)-  v^{(2)}(t) \bigr|_{V_0^*}^2 
	+ \tau \bigl| v^{(1)}(t)-  v^{(2)}(t) \bigr|_{H_0}^2
	+
	\bigl| v^{(1)}_{\Gamma }(t)-  v^{(2)}_{\Gamma }(t) \bigr|_{H_\Gamma }^2 
	} \nonumber \\
	& \quad {} + \int_{0}^{t} 
	\bigl| v^{(1)}(s) -v^{(2)}(s ) \bigr|_{V_0}^2 ds 
	+ 2 \int_{0}^{t} 
	\bigl| \nabla _\Gamma v^{(1)}_{\Gamma }(s ) -\nabla _\Gamma v^{(2)}_{\Gamma }(s ) 
	\bigr|_{H_{\Gamma }^d}^2 ds 
	\nonumber\\
	& \le C 
	\left\{ 
	\bigl| v^{(1)}_0-  v^{(2)}_0 \bigr|_{V_0^*}^2 
	+ \tau \bigl| v^{(1)}_0-  v^{(2)}_0 \bigr|_{H_0}^2 
	+
	\bigl| v^{(1)}_{0\Gamma }-  v^{(2)}_{0\Gamma } \bigr|_{H_\Gamma }^2 
	+ \int_{0}^{T} 
	\bigl| f^{(1)}(s ) -f^{(2)}(s ) \bigr|_{L^2(\Omega )}^2 ds 
	\right. \nonumber\\
	& 
	\quad \left. {} + \int_{0}^{T} 
	\bigl| f^{(1)}_{\Gamma }(s ) -f^{(2)}_{\Gamma }(s ) \bigr|_{H_{\Gamma }}^2 ds 
	\right\} \quad {\it for~all}~t \in [0,T].
	\label{pier3}
\end{align} 

The second result deals with the existence of the solution. 
To the aim, we further assume {\pier that}
\begin{itemize}
 \item[{\pier (A5)}] {\pier there} exist positive constants 
$c_0$, $\varrho >0$ such that 
\begin{gather} 
	|s| \le c_0 \bigl( 1+\widehat{\beta }(r) \bigr) \quad 
	\mbox{for all } r \in \mathbb{R}~
	\mbox{and } s \in \beta (r),
	\label{(A4)-1}
	\\ 
	|s| \le c_0 \bigl( 1+\widehat{\beta }_{\Gamma }(r) \bigr) \quad 
	\mbox{for all } r \in \mathbb{R}~
	\mbox{and } s \in \beta_{\Gamma }(r),
	\label{(A4)-2}
	\\
	\bigl |\beta ^\circ (r) \bigr| \le \varrho \bigl |\beta _{\Gamma }^\circ (r) \bigr |+c_0
	\quad 
	\mbox{for all } r \in \mathbb{R};
	\label{(A4)-3}
\end{gather} 
  \item[{\pier (A6)}] {\pier for the initial data 
$\mbox{\boldmath $ v$}_0=(v_0,v_{\Gamma 0}) \in \mbox{\boldmath $ V$}_0$ the compatibility conditions} 
\begin{equation}
	h_* \le \int_{\Gamma }^{} w_{\Gamma } v_{0\Gamma }d\Gamma \le h^*, \quad 
	\widehat \beta (v_0+m_0) \in L^1(\Omega ), \quad 
	\widehat \beta _{\Gamma }(v_{0\Gamma}+m_0) \in L^1(\Gamma)
	\label{comp.}
\end{equation}
{\pier must hold.}
\end{itemize}
The minimal section $\beta ^\circ $ of $\beta $ is specified by
$\beta ^\circ (r):=\{ r^* \in \beta (r) : |r^*|=\min _{s \in \beta (r)} |s| \}$
and the {\pier same definition applies to} $\beta _{\Gamma }^\circ $. 
{\pier The reader can compare these assumptions with the analogous ones in  \cite[(2.17)--(2.21)]{CF14}}. 

{\pier We have to distinguish between the cases $\tau  >0$ and  $\tau =0$. To this aim, we introduce} the additional regularity assumption for $f$: 
\begin{itemize}
 \item[{\pier (A7)}] $f \in H^1(0,T; L^2(\Omega ))$ or 
$f \in L^2 (0,T; H^1(\Omega ))$. 
\end{itemize} 

\paragraph{Theorem 2.2.} 
{\it Let $\tau >0$. {\pier Then, under the assumptions {\rm (A1)}--{\rm (A6)}, 
there exists a unique solution of {\rm (P)}. Moreover,  if $\tau =0$ and  {\rm (A7)} holds, 
then the problem  {\rm (P)} has a unique solution as well.}}

%%%%% Section 2.3. %%%%%
\subsection{Abstract formulation} 

In this subsection, {\pier an abstract formulation of the problem is given. 
We can write the problem as an evolution inclusion governed by a subdifferential 
operator, with essentially the same approach as in} \cite{CF14, KN96, Kub12}. 

The point of emphasis is that our mass constraint \eqref{(D6)} {\pier reads}
$$
	h_* \le \bigl( \mbox{\boldmath $ w$},\mbox{\boldmath $ v$}(t)
	\bigr)_{\mbox{\scriptsize \boldmath $ H$}_0} \le h^* 
	\quad 
	\mbox{for all }
	t \in [0,T],
$$
with $\mbox{\boldmath $ w$}:=(0,w_\Gamma) \in \mbox{\boldmath $ H$}_0$. 
Then, {\pier by introducing} the convex constraint set 
$$ 
	\mbox{\boldmath $ K$}:=
	\bigl\{ \mbox{\boldmath $ z$} \in \mbox{\boldmath $ V$}_0 \ : \  
	h_* \le ( \mbox{\boldmath $ w$},\mbox{\boldmath $ z$}
	)_{\mbox{\scriptsize \boldmath $ H$}_0} \le h^* 
	\bigr\},
$$
{\pier let  $I_{\mbox{\boldmath \scriptsize $ K$}}:\mbox{\boldmath $ H$}_0 \to [0,+\infty ]$ denote}
the indicator function of $\mbox{\boldmath $ K$}$. 
Now, define {\pier the} proper, lower semicontinuous and convex functional 
$\varphi: \mbox{\boldmath $ H$}_0 \to [0,+\infty ]$ by 
$$
	\varphi (\mbox{\boldmath $ z$}) 
	:= \left\{ 
	\begin{array}{l}
	\displaystyle 
	\frac{1}{2} \int_{\Omega }^{} \bigl | \nabla z \bigr|^2 dx + \int_{\Omega }^{} 
	\widehat{\beta }(z+m_0)dx + \frac{1}{2} \int_{\Gamma }^{} 
	\bigl |\nabla _{\Gamma } z_{\Gamma } \bigr |^2 d\Gamma  
	+ \int_{\Gamma }^{} 
	\widehat{\beta }_{\Gamma } (z_{\Gamma }+m_0)d\Gamma \vspace{2mm}\\
	\hfill
	\mbox{if } 
	\mbox{\boldmath $ z$} \in \mbox{\boldmath $ V$}_0, 
	\widehat{\beta }(z+m_0) \in L^1(\Omega )~
	\mbox{and } \widehat{\beta }_{\Gamma }(z_{\Gamma }+m_0) \in L^1(\Gamma), \vspace{2mm}\\
	+\infty \quad \mbox{otherwise}. 
	\end{array} 
	\right. 
$$
Then, {\pier the problem {\rm (P)} can be stated as the {C}auchy problem for  
an evolution inclusion with a {\takeshi perturbation,} namely}
\begin{gather} 
	\mbox{\boldmath $ A$}_\tau \mbox{\boldmath $ v$}'(t)+\partial (\varphi +I_{\mbox{\boldmath \scriptsize $ K$}})
	\bigl( \mbox{\boldmath $ v$}(t) \bigr)
	\ni P \Bigl( \mbox{\boldmath $ f$}(t) - \mbox{\boldmath $ \Pi $}_0 \bigl (\mbox{\boldmath $ v$}(t) \bigr ) \Bigr) 
	\quad \mbox{in } \mbox{\boldmath $ H$}_0,
	\ \mbox{for a.a.\ } t \in (0,T),
	\label{(E)}
	\\
	\mbox{\boldmath $ v$}(0)=\mbox{\boldmath $ v$}_0 
	\quad \mbox{in } \mbox{\boldmath $ H$}_0,
	\label{(I.C.)}
\end{gather} 
where
$\mbox{\boldmath $ A$}_\tau \mbox{\boldmath $ z$}:=(F^{-1}z+\tau z,z_\Gamma)$ for 
$\tau \ge 0$, 
$P\mbox{\boldmath $ z$}:=(P_0z,z_\Gamma-(1/|\Omega |)\int_{\Omega }^{}zdx)$ {\pier and}
$\mbox{\boldmath $ \Pi $}_0(\mbox{\boldmath $ z$})
:=( \pi (z+m_0),\pi_{\Gamma }(z_\Gamma +m_0))$ 
for all $\mbox{\boldmath $ z$} \in \mbox{\boldmath $ H$}_0$. 

{\pier Hence,  let us recall the paper \cite{CV90}} and express our expectation that {\pier \eqref{(E)}--\eqref{(I.C.)}} 
can be solved by the abstract theory of doubly nonlinear evolution inclusions. All this will be discussed in Section~\ref{Existence}. 
{\pier On the other hand, Theorem~2.2 allows a characterization in terms of regularity of the solution and presence of the {L}agrange multipliers.}

{\pier We aim to point out that analogous remarks were emphasized in \cite{CF14} for an {A}llen--{C}ahn equation
with dynamic boundary conditions and mass constraints; the reader can compare the two problems. In connection with 
\cite{CF14}, we also quote the abstract approach carried out in \cite{FK13}, which however does not comply here with 
the structure of \eqref{(E)}--\eqref{(I.C.)}.}

%%%%% Section 3. %%%%%
\section{Continuous dependence}

In this section, we prove Theorem~2.1. 

\paragraph{Proof of Theorem~2.1.}  {\pier For $i=1,2$ let $(\mbox{\boldmath $ v$}^{(i)}, \mbox{\boldmath $ {\pier\xi} $}^{(i)}, \omega ^{(i)}, \lambda^{(i)} )$ be a solution of {\rm (P)} corresponding to the data  $(f^{(i)}$, $f^{(i)}_{\Gamma }$, $v_0^{(i)}$, $v_{0\Gamma }^{(i)})$.} We {\pier consider} the difference between \eqref{(D1)} written for $v^{(1)}(s)$ of 
$\mbox{\boldmath $ v$}^{(1)}(s)=(v^{(1)}(s),v_\Gamma ^{(1)}(s))$ 
and \eqref{(D1)} written for $v^{(2)}(s)$ of 
$\mbox{\boldmath $ v $}^{(2)}(s)=(v^{(2)}(s),v_\Gamma ^{(2)}(s))$ 
at the time {\pier $s\in (0,T)$}. 
Then, we take 
the inner product with
$v^{(1)}(s)-v^{(2)}(s)$ in $H$. Using the monotonicity of 
$\beta $ and the fact 
$\int_{\Omega }^{} ( v^{(1)}(s)-  v^{(2)}(s)) dx =0$,
we obtain 
\begin{align} 
	\lefteqn{ 
	\frac{1}{2} \frac{d}{ds} \bigl| v^{(1)}(s)-  v^{(2)}(s) \bigr|_{V_0^*}^2 
	+ 
	\frac{\tau }{2} \frac{d}{ds} \bigl| v^{(1)}(s)-  v^{(2)}(s) \bigr|_{H_0}^2 
	} \nonumber \\
	& \quad {} +
	\bigl| v^{(1)}(s ) -v^{(2)}(s ) \bigr|_{V_0}^2 
	- \bigl( \partial _\nu v^{(1)}(s ) - \partial _\nu v^{(2)}(s ), 
	v^{(1)}_\Gamma (s ) -v^{(2)}_\Gamma (s ) \bigr)_{H_\Gamma} 
	\nonumber  \\
	& \le  
	\bigl( f^{(1)}(s ) -f^{(2)}(s ), v^{(1)}(s)-  v^{(2)}(s) \bigr)_{H} 
	\nonumber \\
	& \quad {}
	- \Bigl( \pi \bigl (v^{(1)}(s) +m_0\bigr) -  \pi \bigl (v^{(2)}(s) +m_0\bigr) , 
	v^{(1)}(s)-  v^{(2)}(s) \Bigr)_{\! H}, \label{1st}
\end{align} 
for a.a.\ $s \in (0,T)$. 
Moreover, we take the difference between \eqref{(D3)} written for $v^{(1)}_\Gamma (s)$ of 
and \eqref{(D3)} written for $v^{(2)}_\Gamma (s)$ of 
at the time $t=s$, and take 
the inner product with
$v^{(1)}_\Gamma (s)-v^{(2)}_\Gamma (s)$ in $H_\Gamma${\pier ; hence, we can replace the term 
$$- \bigl( \partial _\nu v^{(1)}(s ) - \partial _\nu v^{(2)}(s ), 
	v^{(1)}_\Gamma (s ) -v^{(2)}_\Gamma (s ) \bigr)_{H_\Gamma} $$
with the corresponding quantity in} \eqref{1st}. 
Then, {\pier by exploiting} the monotonicity of  $\beta_\Gamma$ {\pier and} the {L}ipschitz continuities of 
$\pi$ and $\pi _\Gamma$, we obtain
\begin{align*} 
	\lefteqn{ 
	\frac{d}{ds} \left\{ \bigl| v^{(1)}(s)-  v^{(2)}(s) \bigr|_{V_0^*}^2 
	+ 
	\tau \bigl| v^{(1)}(s)-  v^{(2)}(s) \bigr|_{H_0}^2 
	+ 
	\bigl| v^{(1)}_{\Gamma }(s)-  v^{(2)}_{\Gamma }(s) \bigr|_{H_\Gamma }^2 
	\right\} 
	} \nonumber \\
	& \quad  {} + 
	2 \bigl| v^{(1)}(s ) -v^{(2)}(s ) \bigr|_{V_0}^2 
	+  
	2 \bigl| \nabla _\Gamma v^{(1)}_\Gamma (s ) 
	- \nabla _\Gamma v^{(2)}_\Gamma (s ) \bigr|_{H_\Gamma^d}^2 \\
	& \le  
	\bigl| f^{(1)}(s ) -f^{(2)}(s ) \bigr|_{H}^2 
	+ (1+2L)\bigl| v^{(1)}(s)-  v^{(2)}(s) \bigr|_{H_0}^2 
	+ 
	\bigl| f^{(1)}_{\Gamma }(s ) -f^{(2)}_{\Gamma }(s ) \bigr|_{H_\Gamma }^2 
	\nonumber \\
	& \quad {}
	+ 
	(1+2L_\Gamma )\bigl| v^{(1)}_\Gamma (s)-  v^{(2)}_\Gamma (s) \bigr|_{H_\Gamma}^2, 
\end{align*} 
for a.a.\  $s \in (0,T)$. 
If $\tau >0$, {\pier by applying directly} the {G}ronwall lemma, 
it is straightforward to find a constant $C>0$, 
depending only on $L$, $L_{\Gamma}$ and $T$, such that 
the continuous dependence holds. 
If $\tau=0$, {\pier a known compactness inequality (see, e.g., 
\cite[Thm.~16.4, p.~102]{LM}) 
states that for each $\delta >0$ there exists a positive constant $C_\delta$ such that} 
$$
	|z|_{H_0} \le \delta |z|_{V_0} + C_\delta |z|_{V_0^*} 
	\quad {\rm for~all}~ z \in V_0, 
$$
%(see e.g., \cite[Lemma 8, p.~84]{Sim87}).
Therefore, taking $\delta ^2<1/(2+4L)$  {\pier we have}
\begin{align*} 
	\lefteqn{ 
	(1+2L)\bigl| v^{(1)}(s)-  v^{(2)}(s) \bigr|_{H_0}^2
	} \nonumber \\
	& \le  (1+2L) \left\{ 2 \delta^2 \bigl| v^{(1)}(s)-  v^{(2)}(s) \bigr|_{V_0}^2
	+ 2C_\delta ^2 \bigl| v^{(1)}(s)-  v^{(2)}(s) \bigr|_{V_0^*}^2 \right\} 
	\nonumber  \\
	& \le 
	\bigl| v^{(1)}(s)-  v^{(2)}(s) \bigr|_{V_0}^2
	+ \tilde{C} \bigl| v^{(1)}(s)-  v^{(2)}(s) \bigr|_{V_0^*}^2, 
\end{align*} 
for a.a.\ $s \in (0,T)$  {\pier and some constant $\tilde{C}$ depending only on $L$. 
At this point,  we can analogously apply the {G}ronwall lemma and find a 
constant $C>0$, with the same dependencies as above, such that \eqref{pier3} 
holds. Thus, Theorem~2.1 is completely proved.} \hfill $\Box$

%%%%% Section 4. %%%%%
\section{Existence}
\label{Existence}

This section is devoted to the proof of Theorem~2.2. 
We make use of  {\pier {Y}osida approximations for the} maximal monotone operators $\beta $, $\beta _\Gamma $ 
and of well-known results of this theory (see, \cite{Bar10, Bre73, Ken07}). 
For each $\varepsilon \in (0,1]$, we define 
$\beta _\varepsilon, \beta _{\Gamma,\varepsilon}:\mathbb{R} \to \mathbb{R}$, 
along with the associated resolvent operators $J_\varepsilon, J_{\Gamma,\varepsilon}:\mathbb{R} \to \mathbb{R}$
by 
\begin{gather*}
	\beta _\varepsilon (r)
	:= \frac{1}{\varepsilon } \bigl( r-J_\varepsilon (r) \bigr)
	:=\frac{1}{\varepsilon }\bigl( r-(I+\varepsilon \beta )^{-1} (r)\bigr),
	\\
	\beta _{\Gamma, \varepsilon} (r)
	:= \frac{1}{\varepsilon \varrho} \bigl( r-J_{\Gamma,\varepsilon  }(r) \bigr )
	:=\frac{1}{\varepsilon \varrho}\bigl( r- (I+\varepsilon \varrho \beta _\Gamma )^{-1} (r) \bigr)
	\quad \mbox{for all } r \in \mathbb{R},
\end{gather*}
where $\varrho>0$ is {\pier the same constant as in \eqref{(A4)-3}.} 
Note that the two definitions are not symmetric since in the second it is 
$\varepsilon \varrho$ and not directly $\varepsilon $ to be used as approximation parameter. 
Now, we easily have 
$\beta _\varepsilon(0)=\beta _{\Gamma, \varepsilon}(0)=0$. 
Moreover, the related {M}oreau-{Y}osida regularizations $\widehat{\beta }_\varepsilon, 
\widehat{\beta }_{\Gamma,\varepsilon}$
of $\widehat{\beta }, \widehat{\beta}_{\Gamma }:\mathbb{R} \to \mathbb{R}$ fulfill
\begin{gather*}
	\widehat{\beta }_{\varepsilon }(r)
	:=\inf_{s \in \mathbb{R}}\left\{ \frac{1}{2\varepsilon } |r-s|^2+\widehat{\beta }(s) \right\} 
	= 
	\frac{1}{2\varepsilon } \bigl| r-J_\varepsilon (r) \bigr|^2+\widehat{\beta }(J_\varepsilon r)
	= \int_{0}^{r} \beta _\varepsilon (s)ds,
	\\
	\widehat{\beta }_{\Gamma, \varepsilon }(r)
	:=\inf_{s \in \mathbb{R}}\left\{ \frac{1}{2\varepsilon \varrho } |r-s|^2+\widehat{\beta }_{\pier \Gamma}(s) \right\} 
	= \int_{0}^{r} \beta _{\Gamma,\varepsilon} (s)ds
	\quad \mbox{for all } r\in \mathbb{R}.
\end{gather*}
It is well known that $\beta_\varepsilon$ is {L}ipschitz continuous with {L}ipschitz constant 
$1/\varepsilon $ 
and $\beta_{\Gamma, \varepsilon}$ is also {L}ipschitz continuous with constant 
$1/(\varepsilon \varrho)$. In addition, we have the standard {\fukao properties}
\begin{gather*}
	\bigl |\beta _\varepsilon (r) \bigr | \le \bigl |\beta ^\circ (r) \bigr |, \quad 
	\bigl |\beta _{\Gamma ,\varepsilon }(r) \bigr | \le \bigl |\beta _{\Gamma }^\circ (r) \bigr | 
	\quad \mbox{for all } r \in \mathbb{R},
	\\
	0 \le \widehat{\beta }_\varepsilon (r) \le \widehat{\beta }(r), \quad 
	0 \le \widehat{\beta }_{\Gamma, \varepsilon} (r) \le \widehat{\beta }_{\Gamma }(r)
	\quad \mbox{for all } r \in \mathbb{R}.
\end{gather*}
Here, we note that from {\pier the assumptions \eqref{(A4)-1}, \eqref{(A4)-2} and the above properties} we also obtain
\begin{gather} 
	\bigl |\beta _\varepsilon (r) \bigr | \le c_0 \bigl( 1+\widehat{\beta }_\varepsilon (r) \bigr), 
	\label{(A4)-1e}
	\\
	\bigl |\beta _{\Gamma,\varepsilon }(r) \bigr | \le c_0 \bigl( 1+\widehat{\beta }_{\Gamma, \varepsilon  }(r) \bigr) \quad 
	\mbox{for all } r \in \mathbb{R},
	\label{(A4)-2e}
\end{gather}
{\pier with the same constant $c_0$.} %  as in  \eqref{(A4)-1}, \eqref{(A4)-2}.  
Moreover, thanks to \eqref{(A4)-3} and \cite[Lemma~4.4]{CC13}, the inequality
\begin{equation}
	\bigl |\beta_\varepsilon (r)\bigr | \le \varrho \bigl |\beta _{\Gamma,\varepsilon} (r)\bigr |+c_0
	\quad 
	\mbox{for all } r \in \mathbb{R},
	\label{(A4)-3e}
\end{equation} 
holds for $\beta _\varepsilon $ and $\beta _{\Gamma, \varepsilon}$.

%%%%% Section 4.1. %%%%%
\subsection{Approximation of the problem}

In this subsection, we consider {\pier the approximation of problem} {\rm (P)} in the case when $\tau >0$. 
{\pier The limiting case as $\tau \to 0$ will be discussed later. 
We introduce the following 
{C}auchy problem: for each $\varepsilon \in (0,1]$} find $\mbox{\boldmath $ v$}_\varepsilon $ 
satisfying
\begin{gather} 
	\mbox{\boldmath $ A$}_\tau \mbox{\boldmath $ v$}_\varepsilon '(t)
	+\partial (\varphi _\varepsilon +I_{\mbox{\boldmath \scriptsize $ K$}})
	\bigl( \mbox{\boldmath $ v$}_\varepsilon (t) \bigr) 
	\ni P \Bigl( \mbox{\boldmath $ f$}(t) - \mbox{\boldmath $ \Pi $}_0\bigl( \mbox{\boldmath $ v$}_\varepsilon (t) \bigr) \Bigr) 
	\nonumber\\ {\pier \hskip7cm \mbox{in } \mbox{\boldmath $ H$}_0,
	\ \mbox{for a.a.\ } t \in (0,T),}
	\label{(E)e}
	\\[0.2cm]
	\mbox{\boldmath $ v$}_\varepsilon (0)=\mbox{\boldmath $ v$}_0 
	\quad \mbox{in } \mbox{\boldmath $ H$}_0,
	\label{(I.C.)e}
\end{gather} 
with $\mbox{\boldmath $ v$}_0=( v_0,v_{0\Gamma} ) \in \mbox{\boldmath $ K$}$
satisfying the compatibility conditions \eqref{comp.}. 
Here, $\varphi_\varepsilon:\mbox{\boldmath $ H$}_0 \to [0, +\infty ]$ is 
defined by
$$ 
	\varphi_\varepsilon (\mbox{\boldmath $ z$}) 
	:=
	\left\{ 
	\begin{array}{l} 
	\displaystyle 
	\frac{1}{2} \int_{\Omega }^{} \bigl |\nabla z \bigr |^2 dx 
	+ \int_{\Omega }^{} 
	\widehat{\beta }_\varepsilon ( z +m_0)dx \vspace{2mm}\\
	\displaystyle 
	{} + \frac{1}{2} \int_{\Gamma }^{} \bigl |\nabla _{\Gamma } z_{\Gamma } \bigr |^2 d\Gamma 
	+ \int_{\Gamma }^{} 
	\widehat{\beta }_{\Gamma,\varepsilon } ( z_{\Gamma }+m_0) d\Gamma 
	+ {\pier \frac{\varepsilon }{2} \int_{\Gamma }^{} |z_{\Gamma}|^2d\Gamma } \quad 
	\mbox{if } \mbox{\boldmath $ z$} \in \mbox{\boldmath $ V$}_0, \vspace{2mm}\\
	+\infty \quad \mbox{if } \mbox{\boldmath $ z$} \in \mbox{\boldmath $ H$}_0 
	\setminus \mbox{\boldmath $ V$}_0,
	\end{array} 
	\right. 
$$
moreover, it is understood that  
{\pier $\mbox{\boldmath $ A$}_\tau \mbox{\boldmath $ z$}:=(F^{-1}z+\tau z,z_\Gamma)$,}
$P\mbox{\boldmath $ z$}:=(P_0z,z_\Gamma-(1/|\Omega |)\int_{\Omega }^{}zdx)$ 
and 
$\mbox{\boldmath $ \Pi $}_0(\mbox{\boldmath $ z$})
:=( \pi (z+m_0),\pi_{\Gamma }(z_\Gamma +m_0))$ 
for all $\mbox{\boldmath $ z$}{\pier{}= (z,z_{\Gamma })} \in \mbox{\boldmath $ H$}_0$.

As a remark, {\pier thanks to the {P}oincar\'e--{W}irtiger inequality for functions with $0$ mean value, there is no need to introduce an approximating term like $(\varepsilon /2) \int_{\Omega }^{} |z|^2dx $
in the expression of  $ \varphi _\varepsilon $ above.}
Denote $\partial _* \varphi _\varepsilon $ by the subdifferential 
of $\varphi _\varepsilon :\mbox{\boldmath $ V$}_0 \to [0,+ \infty]$ 
from
$\mbox{\boldmath $ V$}_0$ to $\mbox{\boldmath $ V$}_0^*$.
From \cite[Lemma~3.1]{CF14}, 
we obtain the characterization of $\partial _* \varphi_\varepsilon $ by
\begin{align} 
	&\bigl \langle 
	\partial _* \varphi _\varepsilon (\mbox{\boldmath $ z$}), 
	\bar{\mbox{\boldmath $ z$}}
	\bigr \rangle_{\mbox{\boldmath \scriptsize $ V$}_0^*,\mbox{\boldmath \scriptsize $ V$}_0}
	% \nonumber \\ & 
	{\pier{}=  \bigl (\nabla z,\nabla \bar{z} \bigr )_{L^2(\Omega )^d}+ 
	\bigl (\beta _\varepsilon (z+m_0),\bar{z} \bigr )_{L^2(\Omega )}+
	\bigl (\nabla _{\Gamma } z_{\Gamma },\nabla _{\Gamma } \bar{z}_{\Gamma } \bigr )_{H_{\Gamma}^{d} }}
	\nonumber \\
	& \quad {}
	+
	\bigl (\beta _{\Gamma ,\varepsilon }(z_{\Gamma }+m_0),\bar{z}_{\Gamma } \bigr )_{H_{\Gamma }}
	+\varepsilon (z_\Gamma , \bar{z}_\Gamma )_{H_\Gamma }
	\quad \mbox{for all }{\pier \mbox{\boldmath $ z$}{\pier{}= (z,z_{\Gamma })},\ } 
	 \bar{\mbox{\boldmath $ z$}}=(\bar{z},\bar{z}_{\Gamma }) 
	\in \mbox{\boldmath $ V$}_0.
	\label{charact.}
\end{align} 
Moreover, there exists a positive constant  $C_\varepsilon $ depending on $\varepsilon >0$ such that 
\begin{equation} 
	\bigl|
	\partial _* \varphi _\varepsilon (\mbox{\boldmath $ z$}) 
	\bigr |_{\mbox{\boldmath \scriptsize $ V$}_0^*} 
	\le C_\varepsilon \bigl(1+\varphi _\varepsilon (\mbox{\boldmath $ z$}) \bigr) \quad 
	\mbox{for all } \mbox{\boldmath $ z$} \in \mbox{\boldmath $ V$}_0.
	\label{crit.}
\end{equation}

Now, we recall the fact that the closure $\overline{\mbox{\boldmath $ K$}}$ 
of $\mbox{\boldmath $ K$}$ in 
$\mbox{\boldmath $ H$}_0$ is characterized by 
$$
	\overline{\mbox{\boldmath $ K$}}=
	\bigl\{ \mbox{\boldmath $ z$} \in \mbox{\boldmath $ H $}_0\ : \  
	h_* \le (\mbox{\boldmath $ w$},\mbox{\boldmath $ z$})_{\mbox{\boldmath \scriptsize $ H$}_0}
	\le h^* \bigr\},
$$
which is closed convex subset of $\mbox{\boldmath $ H$}_0$. 
Moreover, there exists {\pier a function} $z_c \in C^1(\overline{\Omega })$ such that 
$$
	\int_{\Omega }^{}z_c dx =0, \quad {z_c}_{|_\Gamma }=\frac{1}{\sigma _0},
$$
{\pier whence} $\mbox{\boldmath $ z$}_c:=(z_c ,1/\sigma _0) \in \mbox{\boldmath $ V$}_0$. 
Then, we {\pier can deduce the following result.}

\paragraph{Proposition 4.1.} 
{\it Let $\tau >0$. 
For each $\varepsilon \in (0,1]$, there 
{\pier exist} a unique 
$$
	\mbox{\boldmath $ v$}_\varepsilon \in H^1(0,T;\mbox{\boldmath $ H$}_0) 
	\cap 
	L^\infty (0,T;\mbox{\boldmath $ V$}_0)
$$
and a pair of functions 
$\mbox{\boldmath $ v$}_\varepsilon^* \in L^2(0,T;\mbox{\boldmath $ H$}_0)$ and 
$\lambda _\varepsilon \in L^2(0,T)$ such that 
$$
	\mbox{\boldmath $ u$}_\varepsilon (t) 
	\in \overline{\mbox{\boldmath $ K$}}
	\quad {\it for~all~} t \in [0,T],
$$ 
and} 
\begin{gather} 
	\mbox{\boldmath $ A$}_\tau \mbox{\boldmath $ v$}_\varepsilon '(t)
	+ \mbox{\boldmath $ v$}_\varepsilon ^* (t) 
	+ \lambda _ \varepsilon (t)\mbox{\boldmath $ w$}
	= P \Bigl( \mbox{\boldmath $ f$}(t) - 
	\mbox{\boldmath $ \Pi $}_0 \bigl( \mbox{\boldmath $ v$}_\varepsilon (t) \bigr) \Bigr)
	\quad \mbox{in } \mbox{\boldmath $ H$}_0, 
	\  \mbox{for a.a.\ } t \in (0,T),
	\label{(L1)}
	\\
	\mbox{\boldmath $ v$}_\varepsilon ^* (t)
	:=\bigl( v_\varepsilon ^* (t),v_{\Gamma, \varepsilon} ^* (t) \bigr) = \partial \varphi _\varepsilon 
	\bigl(\mbox{\boldmath $ v$}_\varepsilon (t) \bigr)
	\quad \mbox{in } \mbox{\boldmath $ H$}_0, 
	\  \mbox{for a.a.\ } t \in (0,T),
	\label{(L2)}
	\\ 
	\lambda _\varepsilon (t)\mbox{\boldmath $ w$}:=
	\lambda _\varepsilon (t) (0,w_\Gamma) \in 
	\partial I_{ \overline{ \mbox{\boldmath \scriptsize $ K$}} }
	\bigl(\mbox{\boldmath $ v$}_\varepsilon (t) \bigr)
	\quad \mbox{in } \mbox{\boldmath $ H$}_0, 
	\  \mbox{for a.a.\ } t \in (0,T),
	\label{(L3)}
	\\
	\mbox{\boldmath $ v$}_\varepsilon (0)=\mbox{\boldmath $ v$}_0 
	\quad \mbox{in } \mbox{\boldmath $ H$}_0.
	\label{(L4)}
\end{gather}

\paragraph{Proof.} We sketch the basic steps of the proof. 
\subparagraph{1.} {\pier We claim that for} a given $\bar{\mbox{\boldmath $ v$}} \in C([0,T];\mbox{\boldmath $ H$}_0)$ 
there exists a unique $${\pier \mbox{\boldmath $ v$} \in  H^1(0,T;\mbox{\boldmath $ H$}_0) 
	\cap 
	L^\infty (0,T;\mbox{\boldmath $ V$}_0) \subset 
	{\takeshi C \bigl ([0,T];\mbox{\boldmath $ H$}_0 \bigr)} }$$ 
 such that
\begin{gather*}
	\mbox{\boldmath $ A$}_\tau \mbox{\boldmath $ v$} '(t)
	+\partial (\varphi _\varepsilon +I_{\mbox{\boldmath \scriptsize $ K$}})
	\bigl( \mbox{\boldmath $ v$} (t) \bigr) 
	\ni P\Bigl( \mbox{\boldmath $ f$}(t)
	-\mbox{\boldmath $ \Pi $}_0 \bigl( \bar{\mbox{\boldmath $ v$} }(t) \bigr)
	\Bigr) 
	\quad \mbox{in } \mbox{\boldmath $ H$}_0,
	\ \mbox{for a.a.\ } t \in (0,T),
\\
	\mbox{\boldmath $ v$}(0)=\mbox{\boldmath $ v$}_0 
	\quad \mbox{in } \mbox{\boldmath $ H$}_0.
\end{gather*}
{\pier Indeed, it suffices to apply} the abstract theory of doubly nonlinear evolution inclusions
(see, e.g., {\pier \cite[Thm.~2.1]{CV90}}). {\pier We point out that, thanks to $\tau >0$, the operator 
$\mbox{\boldmath $ A$}_\tau$ is coercive in  $\mbox{\boldmath $ H$}_0$.}
Then, we construct the map 
$$
	\Psi : \bar{\mbox{\boldmath $ u$}} \mapsto \mbox{\boldmath $ u$},
$$
from $C([0,T];\mbox{\boldmath $ H$}_0)$ into itself.

\subparagraph{2.} For given 
$\bar{\mbox{\boldmath $ u$}}^{(i)} \in C([0,T];\mbox{\boldmath $ H$}_0)$, 
put $\mbox{\boldmath $ u$}^{(i)}:=\Psi \bar{\mbox{\boldmath $ u$}}^{(i)}$ for $i=1,2$. 
{\pier Then, using the monotonicity of $\partial (\varphi _\varepsilon +I_{\mbox{\boldmath \scriptsize $ K$}})$ 
and the special form of $\mbox{\boldmath $ A$}_\tau$, it is not difficult to deduce the estimate} 
\begin{equation}
	\bigl| \mbox{\boldmath $ u$}^{(1)}(t) 
	- \mbox{\boldmath $ u$}^{(2)}(t) \bigr|_{\mbox{\boldmath \scriptsize $ H$}_0}^2
	\le {\pier C_{\tau}}
	\int_{0}^{t} 
	\bigl| \bar{\mbox{\boldmath $ u$}}^{(1)}(s ) 
	-\bar{\mbox{\boldmath $ u$}}^{(2)}(s ) \bigr|_{\mbox{\boldmath \scriptsize $ H$}_0}^2 ds 
	\quad \mbox{for all } t \in [0,T],
\label{pier4}
\end{equation}
{\pier where $  {\pier C_{\tau}}$ is a constant depending on $L, \, L_{\Gamma}$ and $\tau$.
Owing to \eqref{pier4},  we can prove that there exist a} suitable $k \in \mathbb{N}$ such that $\Psi ^k$ is a 
contraction mapping in $C([0,T];\mbox{\boldmath $ H$}_0)$,  Hence, {\pier being $\tau >0$ there 
exists a unique fixed point for $\Psi$ which yields} the unique solution $\mbox{\boldmath 
$ v$}_\varepsilon $ of the  problem \eqref{(E)e}--\eqref{(I.C.)e}. 

\subparagraph{3.} The third step is essentially {\pier the same as in the abstract theory developed in}~\cite{FK13}. 
Put 
$$
	\mbox{\boldmath $ y$}_\varepsilon (t):=
	-\mbox{\boldmath $ A$}_\tau \mbox{\boldmath $ v$}_\varepsilon '(t)
	+P \Bigl( \mbox{\boldmath $ f$}(t)- \mbox{\boldmath $ \Pi $}_0 
	\bigl( \mbox{\boldmath $ v$}_\varepsilon (t) \bigr) \Bigr) 
	\quad \mbox{in } \mbox{\boldmath $ H$}_0,
	\ \mbox{for a.a.\ } t \in (0,T).
$$
{\pier and observe that} $\mbox{\boldmath $ y$}_\varepsilon \in L^2(0,T;\mbox{\boldmath $ H$}_0)$. 
In general, {\pier for each $\mbox{\boldmath $ z$}\in \mbox{\boldmath $ V$}_0$
we have that 
$$\partial (\varphi _\varepsilon +I_{\mbox{\scriptsize \boldmath $ K$}})(\mbox{\boldmath $ z$}) 
\subset \partial _* (\varphi _\varepsilon +I_{\mbox{\scriptsize \boldmath $ K$}})(\mbox{\boldmath $ z$})
= \partial _* \varphi _\varepsilon(\mbox{\boldmath $ z$}) 
+ \partial _*  I_{\mbox{\scriptsize \boldmath $ K$}}(\mbox{\boldmath $ z$}).$$
}%
Thus, there exists 
$\mbox{\boldmath $ v$}_\varepsilon ^{**}(t) \in \partial _* 
I_{\mbox{\scriptsize \boldmath $ K$}} (\mbox{\boldmath $ v$}_\varepsilon (t))$ such that 
$$
	\mbox{\boldmath $ y$}_\varepsilon (t)
	=\partial _* \varphi _\varepsilon 
	\bigl (\mbox{\boldmath $ v$}_\varepsilon (t) \bigr )
	+ \mbox{\boldmath $ v$}_\varepsilon ^{**}(t)
	\quad \mbox{in } \mbox{\boldmath $ V$}_0^*, 
	\  \mbox{for a.a.\ } t \in (0,T).
$$
Moreover, {\pier taking advantage of} \cite[{\pier Prop.~2}]{FK13} {\pier and}
using $\mbox{\boldmath $ z$}_c=(z_c,1/\sigma _0) \in \mbox{\boldmath $ V$}_0$,
we {\pier set}
\begin{equation} 
	\lambda _\varepsilon (t)
	:=(\mbox{\boldmath $ y$}_\varepsilon (t), \mbox{\boldmath $ z$}_c
	)_{\mbox{\scriptsize \boldmath $ H$}_0}
	- \Bigl \langle \partial _* \varphi _\varepsilon 
	\bigl (\mbox{\boldmath $ v$}_\varepsilon (t) \bigr ),\mbox{\boldmath $ z$}_c 
	\Bigr \rangle _{\! \mbox{\scriptsize \boldmath $ V$}_0^*,\mbox{\scriptsize \boldmath $ V$}_0}	
	\quad {\pier \mbox{for a.a.\ } t \in (0,T)}
	\label{(L5)}
\end{equation} 
{\pier and obtain}
$$
	\mbox{\boldmath $ v$}_\varepsilon ^{**}(t)=\lambda _\varepsilon (t) \mbox{\boldmath $ w$} 
	\in \partial I_{ \overline{ \mbox{\scriptsize \boldmath $ K$}} }
	(\mbox{\boldmath $ v$}_\varepsilon (t)) 
	\quad \mbox{in } \mbox{\boldmath $ H$}_0,	
	\  \mbox{for a.a.\ } t \in (0,T),
$$
{\pier where $\mbox{\boldmath $ w$} =(0,w_\Gamma)$ (cf.\ (A4)). Note that 
$\lambda _\varepsilon \in L^2(0,T)$ thanks to \eqref{(L5)} and \eqref{charact.}. As a consequence, both 
$\mbox{\boldmath $ v$}_\varepsilon ^{**}$  and 
$\mbox{\boldmath $ v$}_\varepsilon ^*:=\partial _* \varphi (\mbox{\boldmath $ v$}_\varepsilon )$ are in $  L^2(0,T;\mbox{\boldmath $ H$}_0)$  and  \eqref{(L1)}--\eqref{(L3)} follow with the right regularity.} 
\hfill $\Box$ 
\vspace{2mm}

Let $\tau>0$. Using Proposition~4.1 with 
the characterization \eqref{charact.} of 
$\partial _* \varphi _\varepsilon$
we obtain the following weak formulation: 
\begin{align}
	\lefteqn{ 
	\int_{\Omega }^{} F^{-1} \left( \frac{\partial v_\varepsilon }{\partial t}(t)\right)z dx 
	+ \tau \int_{\Omega }^{} \frac{\partial v_\varepsilon }{\partial t}(t) z dx 
	+ \int_{\Gamma }^{} \frac{\partial v_{\Gamma,\varepsilon } }{\partial t}(t) z_\Gamma d \Gamma 
	+ \int_{\Omega }^{} \nabla v_\varepsilon (t) \cdot \nabla z dx 
	} \nonumber \\
	& \quad {} 
	+ \int_{\Gamma }^{} \nabla_\Gamma  v_{\Gamma,\varepsilon }(t) \cdot \nabla_\Gamma z_\Gamma d\Gamma 
	+ \int_{\Omega }^{} q_\varepsilon (t) z dx 
	+ \int_{\Gamma }^{} q_{\Gamma,\varepsilon }(t) z_\Gamma d\Gamma 
	+ \int_{\Gamma }^{} \lambda _\varepsilon (t)w_\Gamma z_\Gamma d\Gamma 
	\nonumber \\
	& = 0
	\quad \mbox{for all } 
	\mbox{\boldmath $ z$}:=(z, z_\Gamma ) \in \mbox{\boldmath $ V$}_0,
	\label{v.i.e}
\end{align} 
where  
\begin{gather*}
	q_\varepsilon :=\beta _\varepsilon (v_\varepsilon +m_0)
	+ \pi (v_\varepsilon +m_0) - f
	\in L^2(0,T;L^2(\Omega )),
	\\
	q_{\Gamma, \varepsilon} :=\varepsilon v_{\Gamma,\varepsilon }
	+ \beta _{\Gamma,\varepsilon }(v_{\Gamma,\varepsilon }+m_0)
	+ \pi_\Gamma  (v_{\Gamma, \varepsilon}+m_0)-f_\Gamma 
	\in L^2(0,T;H_\Gamma ).
\end{gather*}
{\pier We also introduce the auxiliary quantity}
\begin{equation}
	\omega_\varepsilon (t)
	:= \frac{1}{|\Omega |} 
	\int_{\Omega }^{} q_\varepsilon (t) dx + 
	\frac{1}{|\Omega |} 
	\int_{\Gamma }^{} \left( \frac{\partial v_{\Gamma,\varepsilon } }{\partial t}(t) 
	+ q_{\Gamma,\varepsilon } (t)+\lambda_\varepsilon  (t)w_\Gamma 
	\right) d\Gamma 	
	\label{omega-e}
\end{equation} 
for a.a.\ $t \in (0,T)$. 
By noting that $\partial v_{\Gamma,\varepsilon}/\partial t$  {\pier and $\lambda _\varepsilon w_\Gamma $ lie in $ L^2(0,T;H_\Gamma )$, it turns out that} $\omega _\varepsilon \in L^2(0,T)$. 
Moreover, {\pier according to} \cite[{\pier Prop.~3.2}]{CF14}, for each $\varepsilon \in (0,1]$ {\pier we can infer that}
$v_\varepsilon \in L^2 (0,T;H^2 (\Omega ))$ and 
$v_{\Gamma, \varepsilon} \in L^2 ( 0,T;H^2(\Gamma))$. 
By virtue of this regularity, our approximate problem can be written as 
\begin{gather} 
	F^{-1} \left( \frac{\partial v_\varepsilon }{\partial t} \right) 
	+ \tau \frac{\partial v_\varepsilon }{\partial t} -\Delta v_\varepsilon 
	+ q_\varepsilon  = \omega _\varepsilon 
	\quad \mbox{a.e.\ in } Q,
	\label{(D1,2)e}\\
	v_{\Gamma, \varepsilon }={v_\varepsilon }_{|_\Gamma },
	\quad 
	\partial_\nu v_\varepsilon 
	+ \frac{\partial v_{\Gamma, \varepsilon  }}{\partial t}
	- \Delta _{\Gamma } v_{\Gamma, \varepsilon }
	+ q_{\Gamma, \varepsilon }
	+ \lambda_\varepsilon w _{\Gamma } = 0
	\quad \mbox{a.e.\ on } \Sigma,
	\label{(D3,4)e}\\
	v_\varepsilon (0)=v_0 \quad \mbox{a.e.\ in } \Omega, \quad 
	v_{\Gamma, \varepsilon}(0)=v_{0\Gamma } \quad \mbox{a.e.\ on } \Gamma,
	\label{(D5)e} \\
	h_* \le h_\varepsilon (t):=
	\int_{\Gamma }^{} w_{\Gamma } v_{\Gamma, \varepsilon}(t) d\Gamma 
	\le h^* 
	\quad \mbox{for all }t \in [0,T],
	\label{(D6)e}\\
	\lambda _\varepsilon(t) \in \partial I_{[h_*,h^*]} \bigl (h_\varepsilon (t) \bigr )
	\quad \mbox{for a.a.\ } t \in (0,T).
	\label{(D7)e}
\end{gather} 
{\pier Due to the regularity of the solution, 
$\mbox{\boldmath $ v$}_\varepsilon (t)$ is in $ \overline{ \mbox{\boldmath $ K$}}$
for all $t\in [0,T]$.  Another remark is that the last condition \eqref{(D7)e}} is equivalent to (see, e.g., \cite[Remark~3.2]{CF14})
\begin{equation}
	\lambda _\varepsilon (t) \mbox{\boldmath $ w$} \in 
	\partial I_{\overline{ \mbox{\boldmath \scriptsize $ K$} }} 
	\bigl(\mbox{\boldmath $ v$}_\varepsilon (t) \bigr)
	\quad \mbox{in } \mbox{\boldmath $ H$}_0,
	\ \mbox{for a.a.\ } t \in (0,T).
	\label{(D7)e2}
\end{equation}

%%%%% Section 4.2. %%%%%
\subsection{A priori estimates}

Let $\tau >0$. 
In this subsection, we obtain the uniform estimates independent of $\varepsilon >0$. 
Moreover, our second objective 
{\pier will be to study the limiting behavior as} $\tau \to 0$. 
Therefore, under the additional regularity assumption {\takeshi {\rm (A7)}} 
for $f$ we also obtain {\pier some} uniform estimates independent of 
$\varepsilon >0$ and $\tau >0$.

\paragraph{Lemma 4.1.}
{\it There exist a positive constant $M_1$, independent of $\varepsilon \in (0,1]$, such that
\begin{gather*}
	|v_\varepsilon|_{H^1(0,T;V_0^*)}
	+ \tau ^{1/2}|v_\varepsilon |_{H^1(0,T;H_0)}
	+ |v_\varepsilon |_{L^\infty (0,T;V_0)}
	+ \sup_{t \in (0,T)} \int_{\Omega }^{} \widehat{\beta }_\varepsilon 
	\bigl( v_\varepsilon(t) +m_0 \bigr) dx \\
	{} +
	|v_{\Gamma, \varepsilon}|_{H^1(0,T;H_\Gamma )}+|v_{\Gamma, \varepsilon}|_{L^\infty (0,T;V_\Gamma )}
	+ \sup_{t \in (0,T)} \int_{\Gamma}^{} \widehat{\beta }_{\Gamma, \varepsilon} 
	\bigl( v_{\Gamma, \varepsilon} (t)+m_0 \bigr) d\Gamma \le M_1.
\end{gather*}
Moreover, if {\takeshi {\rm (A7)}} is assumed, then $M_1>0$ is obtained independent of $\varepsilon \in (0,1]$ and $\tau >0$. }

\paragraph{Proof.} {\pier We test \eqref{(D1,2)e} by  $v_\varepsilon '=\partial v_\varepsilon/\partial t \in L^2(0,T;H)$. 
Moreover, we add $v_{\Gamma ,\varepsilon }$ to  both sides of \eqref{(D3,4)e} and use it as the 
boundary condition, obtaining}
\begin{align}
	\lefteqn{ 
	\int_{0}^{t} \bigl |v_\varepsilon '(s ) \bigr |_{V_0^*}^2 ds 
	+ \tau \int_{0}^{t} \bigl| v_\varepsilon '(s) \bigr |_{H_0}^2 ds 
	+ \frac{1}{2} \bigl |v_\varepsilon (t) \bigr |_{V_0}^2 
	+ \int_{\Omega }^{} \widehat{\beta }_\varepsilon 
	\bigl (v_\varepsilon (t)+m_0 \bigr) dx
	} \nonumber \\
	& \quad {} 
	+ \int_{0}^{t} \bigl |v_{\Gamma, \varepsilon }'(s ) \bigr |_{H_\Gamma }^2 ds 
	+ \frac{1}{2} \bigl |v_{\Gamma, \varepsilon }(t) \bigr |_{V_\Gamma }^2 
	+ \int_{\Gamma }^{} \widehat{\beta }_{\Gamma, \varepsilon} 
	\bigl (v_{\Gamma, \varepsilon} (t) +m_0 \bigr) d\Gamma
	+ \frac{\varepsilon }{2} \bigl |v_{\Gamma, \varepsilon }(t) \bigr |_{H_\Gamma}^2 
	\nonumber \\
	& \quad  {}
	+ \int_{0}^{t} \lambda _\varepsilon (s ) 
	\left\{ 
	\int_{\Gamma }^{}
	w_\Gamma v_{\Gamma, \varepsilon} '(s )d\Gamma 
	\right\} ds \nonumber \\
	& \le 
	\frac{1}{2} |v_0|_{V_0}^2 + 
	\int_{\Omega }^{}
	\widehat{\beta }_\varepsilon (v_0+m_0) dx  
	+ \frac{1}{2} |v_{0\Gamma}|_{V_\Gamma }^2 
	+ \int_{\Gamma }^{} \widehat{\beta }_{\Gamma,\varepsilon  }
	(v_{0\Gamma} + m_0) d\Gamma 
	+ \frac{\varepsilon }{2} |v_{0\Gamma} |_{H_\Gamma }^2 
	\nonumber \\
	& \quad  {}
	+ \int_{0}^{t} \Bigl ( {\pier f(s ) %+ v_\varepsilon (s)
	-\pi  \bigl (v_\varepsilon (s ) + m_0
	\bigr )} ,v_\varepsilon '(s ) 
	\Bigr )_{\! H} ds 
	\nonumber \\
	& \quad {}
	+
	\int_{0}^{t} \Bigl ( f_\Gamma (s )
	+v_{\Gamma, \varepsilon}(s )
	-\pi _\Gamma 
	\bigl (v_{\Gamma , \varepsilon} (s ) + m_0
	\bigr ) ,v_{\Gamma, \varepsilon} '(s ) 
	\Bigr )_{\! H_\Gamma } ds \label{3.6-1}
\end{align} 
for all $t \in [0,T]$. 
We note that (cf.\ \eqref{comp.})
\begin{gather}
	\int_{\Omega }^{}
	\widehat{\beta }_\varepsilon (v_0+m_0) dx 
	\le 
	\int_{\Omega }^{}
	\widehat{\beta }(v_0+m_0) dx<+\infty,
	\label{3.6-2}\\
	\int_{\Gamma }^{} \widehat{\beta }_{\Gamma,\varepsilon  }
	(v_{0\Gamma} +m_0) d\Gamma 
	\le
	\int_{\Gamma }^{} \widehat{\beta }_\Gamma
	(v_{0\Gamma} +m_0) d\Gamma <+ \infty.
	\label{3.6-3}
\end{gather}
Also by the chain rule differentiation lemma (see, e.g., 
\cite[{\pier Lemma~4.4, p.~158}]{Bar10} {\pier or} 
\cite[{\pier Lemme~3.3, p.~73}]{Bre73}) 
{\pier and in view of \eqref{(D6)e}--\eqref{(D7)e},} the last term on the left hand side is exactly 
\begin{equation} 
	\int_{0}^{t} \lambda _\varepsilon (s) h_\varepsilon '(s) ds 
	= I_{[h_*,h^*]} \bigl( h_\varepsilon(t) \bigr) -I_{[h_*,h^*]}(h_0) \equiv 0 
	\quad \mbox{for all }  t \in [0,T],
	\label{3.6-4}
\end{equation} 
where $h_0:=(w_\Gamma ,v_{0\Gamma})_{H_\Gamma}$.
We easily see that 
there exists a positive constant $\tilde{M}_1$, depending on 
$L$, $L_\Gamma$, {\pier $|\pi{\pier (m_0)}|$, $|\pi _\Gamma {\pier (m_0)}|$,} $|\Omega |$  
and $|\Gamma |$ {\pier  (but
independent of $\varepsilon \in (0,1]$ and $\tau >0$)}, such that 
\begin{align}
	\lefteqn{  
	\int_{0}^{t} \Bigl ( {\pier f(s )   %+v_\varepsilon (s)
	-\pi  \bigl (v_\varepsilon (s ) + m_0
	\bigr )} ,v_\varepsilon '(s ) 
	\Bigr )_{\! H} ds 
	} \nonumber \\
	& \le \frac{\tau }{2} \int_{0}^{t} \bigl |v_\varepsilon '(s) \bigr |_{H_0}^2 ds 
	+ {\pier \frac{1}{\tau } \int_{0}^{t} 
	\left( \bigl |f(s) \bigr |_{H}^2  
	% + \bigl | v_\varepsilon (s) \bigr |_{H_0}^2 
	+ \Bigl | \pi \bigl (v_\varepsilon (s) + m_0\bigr ) \Bigr |_{H}^2 \right) 
	ds } \nonumber \\
	& \le \frac{\tau }{2} \int_{0}^{t} \bigl |v_\varepsilon '(s) \bigr |_{H_0}^2 ds 
	+ \frac{\tilde{M}_1}{\tau} \int_{0}^{t} 
	\Bigl( 1 +  \bigl |f(s) \bigr |_H^2 
	+ {\pier \bigl | v_\varepsilon (s) \bigr |_{V_0}^2} \Bigr) 
	ds \label{key}
\end{align} 
{\pier and 
\begin{align}
	\lefteqn{ 
	\int_{0}^{t}  \Bigl ( f_\Gamma (s )
	+v_{\Gamma, \varepsilon }(s)-\pi _\Gamma 
	\bigl (v_{\Gamma , \varepsilon} (s ) + m_0
	\bigr ) ,v_{\Gamma, \varepsilon} '(s ) 
	\Bigr )_{\! H_\Gamma } ds 
	} \nonumber \\
	& \le  \frac{1}{2} \int_{0}^{t} \bigl |v_{\Gamma, \varepsilon} '(s) \bigr |_{H_\Gamma}^2 ds 
	+ \tilde{M}_1 \int_{0}^{t} 
	\Bigl( 1 +  \bigl |f_\Gamma (s) \bigr |_{H_\Gamma }^2 
	+ \bigl | v_{\Gamma, \varepsilon }(s) \bigr |_{H_\Gamma}^2 \Bigr) 
	ds  \label{pier6}  % \quad \mbox{for all } t \in [0,T]. 
\end{align} 
for all $t\in [0,T]$.}  {\pier Now, we collect the 
information in \eqref{3.6-2}--\eqref{pier6} and then apply the {G}ronwall lemma
to the inequality resulting from \eqref{3.6-1}. Hence, we prove the lemma in this case and we see 
from \eqref{key} that} the constant $M_1$ depends on $\tau >0$. 

{\pier On the contrary, if} {\takeshi {\rm (A7)}} is assumed, the key estimate \eqref{key} is modified. 
Thanks to the Young inequality, we see that 
\begin{equation} 
	\int_{0}^{t} 
	\Bigl ({\pier - \pi \bigl (v_\varepsilon (s) + m_0 \bigr )},v_\varepsilon '(s ) 
	\Bigr )_{\! H} ds 
	\le 
	\delta \int_{0}^{t} \bigl |v_\varepsilon '(s) \bigr|_{V_0^*}^2 ds
	+ \frac{\tilde{M}_1}{\delta }  \int_{0}^{t} 
	\Bigl( 1 + \bigl | v_\varepsilon (s) \bigr |_{V_0}^2 \Bigr)ds , 
	\label{key2}
\end{equation} 
for all $\delta >0$. 
If we assume $f \in H^1(0,T; L^2(\Omega ))$,  then {\pier we can integrate 
by parts and use the {Y}oung inequality and \eqref{PW}, as follows:}
\begin{align*} 
	\lefteqn{
	\int_{0}^{t} \bigl ( f(s ),v_\varepsilon '(s ) 
	\bigr )_{\!H} ds 
	} \nonumber \\
	& = {}
	- \int_{0}^{t} \bigl ( f'(s ),v_\varepsilon (s ) 
	\bigr )_{\!H} ds
	+ \bigl ( f(t),v_\varepsilon (t) 
	\bigr )_{\!H}
	- \bigl ( f(0),v_0 
	\bigr )_{\!H}
	\\
	& \le 
	{\pier 
	\frac{1}{2} \int_{0}^{t} \bigl| f'(s) \bigr|_{H}^2 ds 
	+
	\frac{C_0}{2} \int_{0}^{t} \bigl| v_\varepsilon (s) \bigr|_{V_0}^2 ds }
	+ 
	\frac{1}{4} \bigl |v_\varepsilon (t) \bigr |_{V_0}^2 
	+ 
	\frac{1}{4} |v_0|_{H_0}^2 
	+
	({\pier C_0}+1)|f|_{C([0,T];L^2(\Omega ))}^2,
\end{align*} 
for all $t \in [0,T]$. 
Thus, taking $\delta <1$ we can {\pier apply the {G}ronwall lemma}
to obtain the estimate with a certain positive constant $M_1$ independent of $\tau >0$. 
{\pier On the other hand, if} we assume 
$f \in L^2(0,T;H^1(\Omega ))$, 
then {\pier we have}
$$
	\int_{0}^{t} \bigl ( f(s) ,v_\varepsilon '(s ) 
	\bigr )_{\!H} ds 
	\le
	\frac{\delta }{2} \int_{0}^{t} \bigl |v_\varepsilon '(s) \bigr|_{V_0^*}^2 ds 
	+ \frac{1}{2\delta }\int_{0}^{t} \bigl |f(s) \bigr |_{H^1(\Omega )}^2 
	ds \quad \mbox{for all } t \in [0,T]. 
$$
Thus, 
by taking $\delta <2/3$, the {G}ronwall inequality works again 
to the conclusion. \hfill $\Box$ 
\vspace{2mm}

Thanks to the growth conditions 
\eqref{(A4)-1}--\eqref{(A4)-2} (see also \eqref{(A4)-1e}--\eqref{(A4)-2e}), we obtain the 
following estimate.

\paragraph{Lemma 4.2.}
{\it There exist a positive constant $M_2$, 
independent of $\varepsilon \in (0,1]$, such that}
$$
	|\lambda _\varepsilon |_{L^2(0,T)} \le M_2.
$$

\paragraph{Proof.} 
From the expression of $\lambda _\varepsilon $, {\pier given by \eqref{(L5)}, 
we infer that} 
\begin{align*} 
	\lambda _\varepsilon (t) 
	& = {} -
	\int_{\Omega }^{}
	\left\{ 
	F^{-1} \left( \frac{\partial v_\varepsilon }{\partial t}(t) \right) 
	+ \tau \frac{\partial v_\varepsilon }{\partial t}(t) 
	+ q_\varepsilon (t) 
	\right\} z_c dx -
	\int_{\Omega }^{} \nabla v_\varepsilon (t) \cdot \nabla z_c dx 
	\nonumber \\
	& \quad  {}-
	\frac{1}{\sigma _0} \int_{\Gamma }^{} 
	\left\{ \frac{\partial v_{\Gamma,\varepsilon  }}{\partial t}(t) 
	+ q_{\Gamma,\varepsilon }(t) 
	\right\} d\Gamma, 
\end{align*} 
for a.a.\ $t \in (0,T)$. Therefore, 
\begin{align*} 
	|\lambda _\varepsilon |_{L^2(0,T)}^2
	& \le 
	6|z_c|_{H_0}^2
	\int_{0}^{T} 
	\left\{ 
	\Bigl|  
	F^{-1} 
	\bigl( v_\varepsilon '(t) \bigr) 
	\Bigr|_{H_0}^2 
	+ \tau ^2
	\bigl|  
	v_\varepsilon '(t)
	\bigr|_{H_0}^2 
	\right\} 
	dt
	+ 6|z_c|_{C(\overline{\Omega })}^2 \int_{0}^{T}  
	\bigl| q_\varepsilon (t) \bigr|_{L^1(\Omega )}^2 dt 
	\nonumber \\
	& \quad {}
	+ 6| z_c |_{V_0}^2 \int_{0}^{T} \bigl| v_\varepsilon (t) \bigr|_{V_0}^2 dt 
	+ \frac{6}{\sigma _0^2} |\Gamma | \int_{0}^{T} 
	\bigl| v_{\Gamma,\varepsilon} '(t) \bigr|_{H_\Gamma } ^2 dt 
	+ \frac{6}{\sigma _0^2} \int_{0}^{T} \bigl| q_{\Gamma ,\varepsilon }(t)
	\bigr|_{L^1({\pier \Gamma} )}^2dt.
\end{align*} 
{\pier By virtue of}  \eqref{(A4)-1e}--\eqref{(A4)-2e}, 
there exists a positive constant $\tilde{M}_2>0$ depending {\pier only} 
on $c_0$, $L$, $L_\Gamma$, {\pier $|\pi{\pier (m_0)}|$ and $|\pi _\Gamma {\pier (m_0)}|$} 
% independent of $\varepsilon \in (0,1]$, 
such that 
\begin{align*} 
\lefteqn{ 
	\bigl| q_\varepsilon (t) \bigr|_{L^1(\Omega )} 
	} \nonumber \\
	& \le 
	\int_{\Omega }^{} c_0
	\Bigl( 1 + \widehat{\beta }_\varepsilon \bigl( v_\varepsilon (t) +m_0\bigr) \Bigr) dx 
	+ \int_{\Omega }^{} 
	\Bigl\{ L \bigl | v_\varepsilon (t) \bigr |  + 
	\bigl |\pi {\pier (m_0)} \bigr| \Bigr\} dx 
	+ \int_{\Omega }^{} \bigl| f(t) \bigr | dx 
	\\
	& \le \tilde{M}_2
	\left\{  1+
	\int_{\Omega }^{}\widehat{\beta }_\varepsilon \bigl( v_\varepsilon (t) +m_0\bigr) dx 
	+ \bigl |v_\varepsilon (t) \bigr|_{L^1(\Omega )} 
	+ \bigl |f(t) \bigr |_{L^1(\Omega )}
	\right\}
\end{align*} 
and 
\begin{align*}  
	\bigl| q_{\Gamma, \varepsilon} (t) \bigr|_{L^1(\Gamma )} 
	& \le 
	\int_{\Gamma}^{} 
	\varepsilon {\pier | v_{\Gamma, \varepsilon }(t)|} d\Gamma 
	 +
	\int_{\Gamma }^{} c_0
	\Bigl( 1 + \widehat{\beta }_{\Gamma, \varepsilon} \bigl( v_{\Gamma, \varepsilon} (t) 
	+m_0\bigr) \Bigr) d\Gamma  
	\\
	& \quad  {}
	+ \int_{\Gamma }^{} \Bigl\{ L \bigl| v_{\Gamma, \varepsilon} (t) \bigr| 
	+ 
	\bigl |\pi_\Gamma  {\pier (m_0)} \bigr | \Bigr\} d\Gamma 
	+ \int_{\Gamma }^{} \bigl |f_\Gamma (t) \bigr | d \Gamma  
	\\
	& \le  \tilde{M}_2
	\left\{ 
	1
	+ \int_{\Gamma }^{} 
	\widehat{\beta }_{\Gamma, \varepsilon} \bigl( v_{\Gamma, \varepsilon} (t) 
	+m_0\bigr) d\Gamma 
	+\bigl |v_{\Gamma, \varepsilon} (t) \bigr|_{L^1(\Gamma )} 
	+ \bigl |f_\Gamma (t) \bigr |_{L^1(\Gamma )}
	\right\}
\end{align*} 
for a.a.\ $t \in (0,T)$.
Therefore, using Lemma~4.1 and taking {\pier into account that
$$ |F^{-1}(v_\varepsilon '(t))|_{H_0}^2 \le {\pier C_0} 
|  F^{-1} ( v_\varepsilon '(t))|_{V_0}^2={\pier C_0}|v_\varepsilon '(t)|_{V_0^*}^2,$$ 
we can find a positive constant $M_2$, independent of $\varepsilon \in (0,1]$, to prove the assertion}. \hfill $\Box$

\paragraph{Lemma 4.3.}
{\it There exist a positive constant $M_3$, 
independent of $\varepsilon \in (0,1]$, such that}
$$
	|\omega _\varepsilon |_{L^2(0,T)} \le M_3.
$$

\paragraph{Proof.} 
From the expression of $\omega _\varepsilon $, {\pier given} by \eqref{omega-e}, 
we have 
\begin{align*}
	|\omega_\varepsilon |_{L^2(0,T)}^2
	& \le 
	\frac{4}{|\Omega |^2} 
	\int_{0}^{T} \bigl |q_\varepsilon (t) \bigr|_{L^1(\Omega )}^2 dt 
	+ 
	\frac{4}{|\Omega |^2} \int_{0}^{T} 
	\bigl| v_{\Gamma,\varepsilon }'(t) \bigr|_{L^1(\Gamma )}^2dt  
	+ 
	\frac{4}{|\Omega |^2} \int_{0}^{T} 
	\bigl| q_{\Gamma,\varepsilon } (t) \bigr|_{L^1(\Gamma )}^2dt 
	\nonumber 
	\\
	& \quad  {}
	+ 
	\frac{4}{|\Omega |^2} |w_\Gamma |_{L^1(\Gamma )}^2 \int_{0}^{T} 
	\bigl | \lambda_\varepsilon  (t) \bigr|^2 dt.
\end{align*} 
Thus, Lemmas~4.1 and~4.2 {\pier ensure the existence of} a positive constant $M_3${\pier , 
independent of $\varepsilon \in (0,1]$, which yields a bound for 
$|\omega _\varepsilon|_{L^2(0,T)}$}. \hfill $\Box$

\paragraph{Lemma 4.4.}
{\it There exist two positive constants $M_4$ and $M_5$, 
independent of $\varepsilon \in (0,1]$, such that}
\begin{gather*}
	\bigl| \beta _\varepsilon (v_\varepsilon +m_0 ) \bigr |_{L^2(0,T;L^2(\Omega ))} 
	+ \bigl |\beta _\varepsilon (v_{\Gamma,\varepsilon }+m_0) \bigr |_{L^2(0,T;H_\Gamma )} 
	\le M_4,
	\\
	|v_\varepsilon |_{L^2(0,T;H^{3/2}(\Omega ))}
	+ 
	|\partial _\nu v_\varepsilon |_{L^2(0,T;H_\Gamma )} 
	\le M_5.
\end{gather*}

\paragraph{Proof.} Testing \eqref{(D1,2)e} by 
$\beta _\varepsilon (v_\varepsilon +m_0) \in L^2(0,T;H^1(\Omega))$ 
and using \eqref{(D3,4)e}. 
Then, integrating it over $\Omega \times (0,t)$ with respect to $(x,s)$, we infer that 
\begin{align*}
	\lefteqn{ 
	\int_{0}^{t} \!\!\! \int_{\Omega }^{} \beta '_\varepsilon 
	\bigl( v_\varepsilon (s) +m_0 \bigr)
	\bigl | \nabla v_\varepsilon (s) \bigr |^2 dxds
	+ 
	\int_{0}^{t} 
	\Bigl | \beta _\varepsilon 
	\bigl (v_\varepsilon (s) +m_0 \bigr) \Bigr |^2_{L^2(\Omega )} ds 
	} \nonumber \\
	& \quad {} 
	+ \int_{0}^{t} \!\!\! \int_{\Gamma }^{} \beta '_\varepsilon 
	\bigl( v_{\Gamma, \varepsilon} (s) +m_0 \bigr)
	\bigl | \nabla_\Gamma  v_{\Gamma,\varepsilon} (s) \bigr |^2 d\Gamma ds
	\nonumber \\
	& \quad  {}
	+ 
	\int_{0}^{t} \!\!\! \int_{\Gamma }^{} 
	\beta _{\Gamma, \varepsilon } \bigl( v_{\Gamma, \varepsilon }(s) + m_0 \bigr)
	\beta _\varepsilon \big( v_{\Gamma, \varepsilon } (s) +m_0 \bigr)
	d\Gamma ds
	\nonumber \\
	& \le  \int_{0}^{t} \Bigl ( f(s ) - F^{-1} \bigl (v_\varepsilon '(s) \bigr)
	- \tau v_\varepsilon '(s)
	- \pi \bigl (v_\varepsilon (s ) +m_0 
	\bigr ) + \omega_\varepsilon(s) ,\beta _\varepsilon 
	\bigl( v_\varepsilon (s ) +m_0 \bigr) 
	\Bigr )_{\! L^2(\Omega )} ds 
	\\
	& \quad  {}
	+
	\int_{0}^{t} \Bigl ( f_\Gamma (s )-v_{\Gamma ,\varepsilon }'(s)
	-\pi _\Gamma 
	\bigl (v_{\Gamma , \varepsilon} (s ) + m_0
	\bigr ) 
	-\lambda _\varepsilon(s) w_\Gamma ,
	\beta _\varepsilon 
	\bigl( v_{\Gamma, \varepsilon}(s ) + m_0 
	\bigr) 
	\Bigr )_{\! H_\Gamma } ds
	\\
	& \quad  {}
	-\varepsilon
	\int_{0}^{t} \Bigl ( 
	 v_{\Gamma ,\varepsilon } (s),
	\beta _\varepsilon 
	\bigl( v_{\Gamma, \varepsilon}(s ) + m_0 
	\bigr) 
	\Bigr )_{\! H_\Gamma } ds
	\quad \mbox{for all } t\in [0,T],
\end{align*} 
where we should take care that 
$(\beta _\varepsilon (v_\varepsilon +m_0))_{|_\Gamma }
=\beta _\varepsilon (v_{\Gamma, \varepsilon} +m_0) \in L^2(0,T;H^1(\Gamma))$.
Here, we use the assumption \eqref{(A4)-3e} to deduce that 
\begin{align*} 
\lefteqn{ 
	\int_{0}^{t} \!\!\! \int_{\Gamma }^{} 
	\beta _{\Gamma, \varepsilon } \bigl( v_{\Gamma, \varepsilon }(s) +m_0 \bigr)
	\beta _\varepsilon \big( v_{\Gamma, \varepsilon } (s)+m_0 \bigr)
	d\Gamma ds 
	} \nonumber \\
	& = 
	\int_{0}^{t} \!\!\! \int_{\Gamma }^{} 
	\Bigl|
	\beta _{\Gamma, \varepsilon } \bigl( v_{\Gamma, \varepsilon }(s) + m_0 \bigr)
	\Bigr| 
	\Bigl| 
	\beta _\varepsilon \big( v_{\Gamma, \varepsilon } (s) + m_0 \bigr) 
	\Bigr| 
	d\Gamma ds 
	\\
	& \ge
	\frac{1}{\varrho } \int_{0}^{t} \!\!\! \int_{\Gamma }^{} 
	\Bigl | \beta _\varepsilon \bigr( v_{\Gamma ,\varepsilon } (s) + m_0 \bigr) \Bigr | ^2 
	d\Gamma ds
	- \frac{c_0}{\varrho }\int_{0}^{t} \!\!\! \int_{\Gamma }^{} 
	\Bigl |\beta _\varepsilon \bigr( v_{\Gamma ,\varepsilon } (s) +m_0 \bigr) \Bigr| 
	d\Gamma ds 
	\\
	& \ge
	\frac{1}{2\varrho } \int_{0}^{t} 
	\Bigl | \beta _\varepsilon \bigr( v_{\Gamma ,\varepsilon } (s) + m_0 \bigr) \Bigr |_{H_\Gamma } ^2 ds
	- \frac{c_0^2}{2\varrho } T |\Gamma | 
	\quad \mbox{for all } t \in [0,T],
\end{align*} 
because $\beta _\varepsilon (r)$ and $\beta _{\Gamma, \varepsilon }(r)$ have the same sign 
for all $r \in \mathbb{R}$. 
We also note that 
\begin{gather*}
	\int_{0}^{t} \!\!\! \int_{\Omega }^{} \beta '_\varepsilon 
	\bigl( v_\varepsilon (s) + m_0 \bigr)
	\bigl | \nabla v_\varepsilon (s) \bigr |^2 dxds \ge 0,
	\\
	\int_{0}^{t} \!\!\! \int_{\Gamma }^{} \beta '_\varepsilon 
	\bigl( v_{\Gamma, \varepsilon} (s) +m_0 \bigr)
	\bigl | \nabla_\Gamma v_{\Gamma,\varepsilon} (s) \bigr |^2 d\Gamma ds \ge 0
%	\quad \mbox{for all } t\in [0,T].
\end{gather*}
{\pier for all $t\in [0,T].$}
Moreover, using the {Y}oung inequality and the fact $\varepsilon \le 1$ we have
\begin{align*} 
	\lefteqn{ 
	{}- \varepsilon \int_{0}^{t} \Bigl( v_{\Gamma, \varepsilon} (s), 
	\beta _\varepsilon \bigl( v_{\Gamma, \varepsilon} (s) +m_0  \bigr) \Bigr)_{\! H_\Gamma}  ds 
	} \nonumber \\
	& \le \frac{\delta }{2} \int_{0}^{t}   
	\Bigl| \beta _\varepsilon \bigl( v_{\Gamma, \varepsilon} (s) +m_0 \bigr) \Bigr|_{H_\Gamma }^2  ds
	+ \frac{1}{2\delta } \int_{0}^{t} \bigl |v_{\Gamma ,\varepsilon }(s) \bigr |_{H_\Gamma }^2 d\Gamma 
	%\quad \mbox{for all } t\in [0,T],
\end{align*} 
{\pier for all $t \in [0,T]$ and $\delta >0$.} 
Now, there exists a positive constant $\tilde{M}_4$, {\pier which} depends only on ${\pier C_0}$, $L$, $L_\Gamma$, {\pier $|\pi {\pier (m_0)}|$, $|\pi _\Gamma {\pier (m_0)}|$, $|\Omega |$, 
$|\Gamma |$} and $T$, such that
\begin{align*} 
	\lefteqn{ 
	\int_{0}^{t} \Bigl ( f(s ) - F^{-1} \bigl (v_\varepsilon '(s) \bigr)
	- \tau v_\varepsilon '(s)
	- \pi \bigl (v_\varepsilon (s ) +m_0 
	\bigr ) + \omega_\varepsilon(s) ,\beta _\varepsilon 
	\bigl( v_\varepsilon (s ) +m_0 \bigr) 
	\Bigr )_{\! L^2(\Omega )} ds 
	} \nonumber \\
	& \le  
	\frac{1}{2} \int_{0}^{t} 
	\Bigl| \beta _\varepsilon 
	\bigl( v_\varepsilon (s ) +m_0 \bigr) 
	\Bigr|_{L^2(\Omega )}^2 ds
	\\
	& \quad {}+
	\tilde{M}_4
	\left( 
	1+|f|_{L^2(0,T;L^2(\Omega ))}^2 
	+ |v_\varepsilon' |_{L^2(0,T;V_0^*)}^2 
	+ \tau ^2 |v_\varepsilon' |_{L^2(0,T;H_0)}^2 
	+ |v_\varepsilon |_{L^2 (0,T;H_0)}^2 
	+ |\omega _\varepsilon |_{L^2(0,T)}^2 
	\right),
\end{align*} 
and 
\begin{align*} 
	\lefteqn{ 
	\int_{0}^{t} \Bigl ( f_\Gamma (s )-v_{\Gamma ,\varepsilon }'(s)
	-\pi _\Gamma 
	\bigl (v_{\Gamma , \varepsilon} (s ) + m_0
	\bigr ) 
	-\lambda _\varepsilon(s) w_\Gamma ,
	\beta _\varepsilon 
	\bigl( v_{\Gamma, \varepsilon}(s ) + m_0 
	\bigr) 
	\Bigr )_{\! H_\Gamma } ds
	} \nonumber \\
	& \le 
	\frac{\delta }{2} \int_{0}^{t} 
	\Bigl| \beta _\varepsilon 
	\bigl( v_{\Gamma, \varepsilon }(s ) +m_0 \bigr) 
	\Bigr|_{H_\Gamma }^2 ds \\
	& \quad {}
	+ 
	\frac{\tilde{M}_4}{2\delta }
	\left( 
	1+|f_\Gamma |_{L^2(0,T;H_\Gamma )}^2 
	+ |v_{\Gamma,\varepsilon}' |_{L^2(0,T;H_\Gamma )}^2 
	+ |v_{\Gamma,\varepsilon}|_{L^2(0,T;H_\Gamma )}^2 
	 + |\lambda _\varepsilon |_{L^2(0,T)}^2 |w_\Gamma |_{H_\Gamma }^2
	\right),
\end{align*} 
for all $t \in [0,T]$ and $\delta >0${\pier ,  with the help of} the {Y}oung inequality. 
Thus, {\pier choosing} $\delta <1/(2\varrho)$ {\pier and recalling Lemmas~4.1--4.3} we deduce that 
there exist a positive constant $M_4${\pier ,
% depending on $\varrho $, $c_0$, $T$, $|\Gamma|$, 
% $|f|_{L^2(0,T;L^2(\Omega ))}$, 
% $|f_\Gamma|_{L^2(0,T;H_\Gamma )}$, $M_1$, $M_2$, $M_3$ and $\tilde{M}_4$, 
independent of $\varepsilon \in (0,1]$,} such that 
$$
	\bigl |\beta _\varepsilon (v_\varepsilon +m_0)\bigr |_{L^2(0,T;L^2(\Omega ))} 
	+ \bigl |\beta _\varepsilon (v_{\Gamma,\varepsilon }+m_0)\bigr |_{L^2(0,T;H_\Gamma )} 
	\le M_4.
$$
{\pier Next}, we can compare the terms in \eqref{(D1,2)e} and conclude that 
$$
	|\Delta v_\varepsilon |_{L^2(0,T;L^2(\Omega ))}  \ {\pier \hbox{ is bounded independently of } \varepsilon,}
$$
whence, {\pier taking Lemma~4.1 into account} and applying the theory of the elliptic regularity 
(see, e.g., \cite[{\pier Thm.~3.2, p.~1.79}]{BG87}), we have that 
$$
	|v_\varepsilon |_{L^2(0,T;H^{3/2}(\Omega ))}
	\le \tilde M_5,
$$
and, owing to the trace theory (see, e.g., \cite[{\pier Thm.~2.25, p.~1.62}]{BG87}), that
$$
	|\partial _\nu v_\varepsilon |_{L^2(0,T;H_\Gamma )} 
	\le \tilde M_5.
$$
{\pier for some  constant $\tilde M_5$ independent of $\varepsilon \in (0,1]$. }\hfill $\Box$

\paragraph{Lemma 4.5.}
{\it There exist positive constants $M_6$, $M_7$ and $M_8$,
independent of $\varepsilon \in (0,1]$, such that} 
\begin{gather*}
	\bigl |\beta _{\Gamma, \varepsilon} (v_{\Gamma,\varepsilon }+m_0) 
	\bigr |_{L^2(0,T;H_\Gamma )} \le M_6,
	\quad
	|v_{\Gamma, \varepsilon} |_{L^2(0,T;H^2(\Gamma))}
	\le M_7,
	\quad
	|v_\varepsilon |_{L^2(0,T;H^2(\Omega ))}
	\le M_8.
\end{gather*}

\paragraph{Proof.} We test \eqref{(D3,4)e} by 
$\beta _{\Gamma, \varepsilon} (v_{\Gamma, \varepsilon}+m_0) \in L^2(0,T;V_\Gamma)$ and 
integrate on the boundary, deducing that 
\begin{align}
	\lefteqn{ 
	\int_{\Gamma}^{} \widehat{\beta }_{\Gamma, \varepsilon}
	\bigl( v_{\Gamma, \varepsilon} (t) +m_0 \bigr) d\Gamma 
	+ \int_{0}^{t} \!\!\! \int_{\Gamma}^{} \beta '_{\Gamma, \varepsilon }
	\bigl( v_{\Gamma, \varepsilon} (s) +m_0 \bigr)
	\bigl | \nabla_\Gamma v_{\Gamma, \varepsilon} (s) \bigr |^2 d\Gamma ds
	} \nonumber \\
	& \quad {}
	+ \int_{0}^{t}\Bigl| \beta _{\Gamma, \varepsilon }
	\bigl( v_{\Gamma, \varepsilon}(s) +m_0 
	\bigr) \Bigr|_{H_\Gamma }^2  ds \nonumber \\
	& \le 
	\int_{\Gamma}^{} \widehat{\beta }_{\Gamma, \varepsilon}
	\bigl( v_{0\Gamma} +m_0 \bigr) d\Gamma 
	-  \int_{0}^{t} \Bigl({\pier \varepsilon \, v_{\Gamma, \varepsilon } (s)  + \partial_\nu v_{\Gamma , \varepsilon}(s)}   ,  
	\beta _{\Gamma, \varepsilon } \bigl (v_{\Gamma, \varepsilon }(s) +m_0 \bigr) 
	\Bigr)_{\! H_\Gamma } ds \nonumber \\
	& \quad {}
	+
	\int_{0}^{t} 
	\Bigl ( f_\Gamma (s )
	-\pi _\Gamma 
	\bigl (v_{\Gamma , \varepsilon} (s )  {\pier{} +m_0}
	\bigr ) 
	-\lambda _\varepsilon (s) w_\Gamma ,
	\beta _{\Gamma, \varepsilon }
	\bigl( v_{\Gamma, \varepsilon}(s ) {\pier{} +m_0}
	\bigr) 
	\Bigr )_{\! H_\Gamma } ds, \label{last}
\end{align} 
for all $t \in [0,T]$. 
We note that 
$$
	\int_{0}^{t} \!\!\! \int_{\Gamma}^{} \beta '_{\Gamma, \varepsilon }
	\bigl( v_{\Gamma, \varepsilon} (s) +m_0 \bigr)
	\bigl | \nabla_\Gamma v_{\Gamma, \varepsilon} (s) \bigr |^2 d\Gamma ds 
	\ge 0, 
$$
{\pier due to the properties of $\beta _{\Gamma,\varepsilon }$, and}

$$
	\int_{\Gamma }^{} \widehat{\beta }_{\Gamma,\varepsilon }(v_{0\Gamma }+m_0) d\Gamma 
	\le 
	\int_{\Gamma }^{} \widehat{\beta }_{\Gamma}(v_{0\Gamma }+m_0) d\Gamma < +\infty, 
$$
by virtue of \eqref{comp.}{\pier . By} applying the {Y}oung inequality in the 
last two terms {\pier of~\eqref{last},} 
we see that there exist a positive constant $\tilde{M}_6$ 
% depends on 
%$|f|_{L^2(0,T;H_\Gamma )}$, $L_\Gamma $, 
% $|\pi _\Gamma {\pier (m_0)}|$, $|\Gamma |$, $T$, $M_1$, $M_2$, $M_3$, $M_4$ and $M_5$
{\pier independent of $\varepsilon \in (0,1]$} such that 
$$
	\bigl |\beta _{\Gamma, \varepsilon} (v_{\Gamma,\varepsilon }+m_0) 
	\bigr |_{L^2(0,T;H_\Gamma )} \le \tilde{M}_6.
$$
Hence, by comparison in \eqref{(D3,4)e} we also infer that
$$
	|\Delta_\Gamma v_{\Gamma, \varepsilon} |_{L^2(0,T;H_\Gamma )} 
	 \le {\pier \tilde{M}_7}
$$
and consequently (see, e.g., \cite[Section~4.2]{Gri09})
\begin{eqnarray*} 
	|v_{\Gamma, \varepsilon }|_{L^2(0,T;H^2(\Gamma ))} 
	& \le & \left( |v_{\Gamma, \varepsilon }|_{L^2(0,T;V_\Gamma)}^2+
	|\Delta _\Gamma v_{\Gamma, \varepsilon }|_{L^2(0,T;H_\Gamma)}^2 \right)^\frac{1}{2} \\
	& \le & \bigl( M_1^2T+{\pier \tilde{M}_7^2} \bigr) ^\frac{1}{2}  =:  M_7. 
\end{eqnarray*} 
{\pier Then,} in view of Lemma~4.4, using the theory of the 
elliptic regularity (see, e.g., \cite[{\pier Thm.~3.2, p.~1.79}]{BG87} along with the {\pier boundedness of
$|v_{\Gamma, \varepsilon }|_{L^2(0,T;H^{3/2}(\Gamma ))}$,} it turns out that 
$$
	|v_\varepsilon |_{L^2(0,T;H^2(\Omega ))} \le M_8
$$
for some positive constant $M_8$ independent of $\varepsilon \in (0,1]$. 
\hfill $\Box$

\paragraph{Remark 4.1.} All constants $M_k${\pier , for $k$ from $1$ to $8$}, are obtained 
independent{\pier ly} of $\tau >0 $ {\pier provided that}
{\takeshi {\rm (A7)}} is assumed. Actually, under the additional assumption 
{\takeshi {\rm (A7)}} 
the positive constant {\pier $M_1$ in Lemma~4.1 is} {\fukao independent} of 
$\tau >0$. 

%%%%% Section 4.3. %%%%%
\subsection{Passage to the limit as $\varepsilon \to 0$}

In this subsection, 
we keep $\tau >0$ fixed and conclude the existence proof by passage to the limit of the 
approximate solutions as $\varepsilon \to 0$.
Indeed, owing to the {\pier uniform estimates stated in} Lemmas from~4.1 to~4.5, there 
{\pier exist} a subsequence of $\varepsilon $ (not relabeled) and some 
limit functions $v$, $v_\Gamma$, ${\pier\xi}$, ${\pier\xi} _\Gamma$, $\omega$, $\lambda $ such that 
\begin{gather} 
	v_\varepsilon \to v \quad \mbox{weakly star in } 
	H^1(0,T;H_0) \cap L^\infty (0,T;V_0) \cap L^2 \bigl (0,T;H^2(\Omega ) \bigr), 
	\label{wcv}
	\\
	v_{\Gamma, \varepsilon} \to v_\Gamma  \quad \mbox{weakly star in } 
	H^1(0,T;H_\Gamma ) \cap L^\infty (0,T;V_\Gamma ) \cap L^2 \bigl (0,T;H^2(\Gamma ) \bigr),
	\label{wcvg}
	\\
	\beta _\varepsilon (v_\varepsilon+m_0) \to {\pier\xi} \quad \mbox{weakly in } 
	L^2 \bigl( 0,T;L^2(\Omega ) \bigr),
	\label{wcb}
	\\ 
	\beta _{\Gamma ,\varepsilon } (v_{\Gamma, \varepsilon}+m_0) \to {\pier\xi} _\Gamma 
	\quad \mbox{weakly in } 
	L^2(0,T;H_\Gamma ),
	\label{wcbg}
	\\
	\omega _\varepsilon \to \omega  \quad \mbox{weakly in } 
	L^2(0,T),
	\label{wco}
	\\
	\lambda _\varepsilon \to \lambda  \quad \mbox{weakly in } 
	L^2(0,T),
	\label{wcl}
\end{gather} 
as $\varepsilon \to 0$. 
From \eqref{wcv} and \eqref{wcvg}, due to strong compactness results 
(see, e.g., \cite[{\pier Sect.~8, Cor.~4}]{Sim87}) 
we {\pier have that}
\begin{gather} 
	v_\varepsilon \to v \quad \mbox{strongly in } 
	C\bigl( [0,T];H_0 \bigr) \cap L^2 (0,T;V_0),
	\label{scv}
	\\
	v_{\Gamma, \varepsilon} \to v_\Gamma  \quad \mbox{strongly in } 
	C\bigl( [0,T];H_\Gamma \bigr) \cap L^2 (0,T;V_\Gamma ),
	\label{scvg}
\end{gather}
as $\varepsilon \to 0$.
Moreover, on account of \eqref{(D6)e} {\pier and \eqref{wcvg}}
it is a standard matter to deduce that 
\begin{equation} 
	h_\varepsilon \to h \quad \mbox{weakly in } 
	H^1(0,T)
	\ \mbox{and strongly in } C \bigl( [0,T] \bigr),
	\label{wsch}
\end{equation} 
% as $\varepsilon \to 0$, 
{\pier where}
$$
	h_* \le h(t):=\int_{\Gamma }^{} w_\Gamma v_\Gamma (t)d\Gamma \le h^*
	\quad \mbox{for all } t\in [0,T].
$$
We point out that \eqref{(D3,4)e}, \eqref{wcv} and \eqref{wcvg} imply that 
$v_\Gamma =v _{|_\Gamma}$ a.e.\ on $\Sigma$, while \eqref{(D5)e}, \eqref{scv}{\pier , \eqref{scvg}  entail}
$$
	v(0)=v_0 \quad \mbox{a.e.\ in } \Omega, 
	\quad 
	v_\Gamma (0) = v_{0\Gamma } 
	\quad \mbox{a.e.\ on } \Gamma.
$$
Now, \eqref{wcl} and \eqref{wsch} and the maximal monotonicity of $\partial I_{[h_*,h^*]}$
allow us to conclude that 
$$
	\lambda \in \partial I_{[h_*,h^*]}(h)
	\quad \mbox{a.e.\ in } (0,T),
$$
{\pier that is equivalent to \eqref{pier5}.}
Moreover, {\pier \eqref{scv}--\eqref{scvg}} and 
the {L}ipschitz continuity of {\pier $\pi ,\, \pi _\Gamma $} imply that 
\begin{gather*}
	\pi (v_\varepsilon+ m_0 ) \to \pi (v+m_0) \quad \mbox{strongly in }
	C \bigl( [0,T];L^2(\Omega ) \bigr),
	\\
	\pi _\Gamma (v_{\Gamma, \varepsilon}+m_0) \to \pi _\Gamma (v_\Gamma+m_0) \quad \mbox{strongly in }
	C \bigl( [0,T];H_\Gamma \bigr),
\end{gather*} 
as $\varepsilon \to 0$. 
At this point, we can pass to the limit in \eqref{(D1,2)e} and \eqref{(D3,4)e} obtaining {\pier \eqref{(D1)} and \eqref{(D3)}.}
Moreover, by applying \cite[{\pier Prop.~2.2, p.~38}]{Bar10} and using 
\eqref{wcb}--\eqref{wcbg} with \eqref{scv}--\eqref{scvg}, we obtain 
$$
	{\pier\xi} \in \beta (v+m_0)
	\quad \mbox{a.e.\ in } Q, \quad 
	{\pier\xi}_\Gamma \in \beta_{\Gamma } (v_{\Gamma }+m_0)
	\quad \mbox{a.e.\ on } \Sigma.
$$
Thus, it turns out that the pair $\mbox{\boldmath $ v$}=(v,v_\Gamma )$ 
{\pier yields, along with  $\mbox{\boldmath $ \xi$}=(\xi   , \xi_\Gamma )$, $\omega$ and $\lambda $,}
a  solution of the limit problem, which can be stated exactly 
as in \eqref{(D1)}--{\pier\eqref{pier5}}. Also, we 
note the regularities 
$v \in C([0,T];V_0)$ and $u_\Gamma \in C([0,T];V_\Gamma )$ for the 
solution as a consequence of \eqref{wcv}--\eqref{wcvg}.

%%%%% Section 4.4. %%%%%
\subsection{Passage to the limit as $\tau \to 0$}

In this subsection, we discuss the limiting problem as $\tau \to 0$.
We need to assume the additional regularity {\takeshi {\rm (A7)}} for $f$. 
For each {\pier $\tau >0$}, let now $\mbox{\boldmath $ v$}_\tau :=(v_\tau ,v_{\Gamma ,\tau })$ be the 
solution to \eqref{(D1)}--{\pier\eqref{pier5}} with 
related $\omega _\tau$, $\lambda _\tau$ and 
$$
	h_\tau (t):=\int_{\Gamma }^{} w_\Gamma v_{\Gamma,\tau}(t)d\Gamma 
	\quad \mbox{for all } t \in [0,T].
$$
{\pier On account} of Lemma~4.1 with Remark~4.1, we use the uniform estimates {\pier in Lemmas 4.1--4.5
to perform the limit procedure as} $\tau \to 0$. 

{\pier As in the previous passage to the limit as $\varepsilon \to 0$,  also in this case 
a subsequence of $\tau $ (not relabeled) and some 
limit functions $v$, $v_\Gamma$, ${\pier\xi}$, ${\pier\xi} _\Gamma$, $\omega$, $\lambda $ 
can be found in order that the same convergences as in \eqref{wcvg}--\eqref{wcl} and}
\begin{equation} 
	v_\tau \to v \quad \mbox{weakly star in } 
	H^1(0,T;V_0^*) \cap L^\infty (0,T;V_0) \cap L^2 \bigl (0,T;H^2(\Omega ) \bigr)
	\label{wcv2}
\end{equation} 
{\pier hold as $\tau \to 0$.  We can still deduce the same strong convergences as in \eqref{scv}--\eqref{wsch} 
and the passage to the limit can be carried out in a similar way.  
Of course, here}  we have to point out that {\pier (cf.\ the estimate in Lemma~4.1)}
$$
	\tau v_\tau ' \to 0 \quad \mbox{strongly in } L^2(0,T;H_0)
$$
as $\tau \to 0$, which is important when we pass to the limit in the equation 
\eqref{(D1)}, obtaining 
\begin{equation}
	F^{-1}\left( \frac{\partial v}{\partial t} \right) 
	-\Delta v + {\pier\xi} + \pi(v+m_0)= f + \omega 
	\quad \mbox{a.e.\ in } Q,
\label{pier7}
\end{equation}
{\pier to be coupled} with \eqref{(D2)}--{\pier\eqref{pier5}}. 

\paragraph{Remark 4.2.}
{\pier On the side of} the proof, one can make the remark that the solution  
{\pier component $\mbox{\boldmath $ v$}=(v,v_\Gamma )$ of the problem
solves the abstract formulation (see Subsections~2.4 and 4.1)}
\begin{gather*}
	\mbox{\boldmath $ v $} \in H^1(0,T;\mbox{\boldmath $ V $}^*_0) 
	\cap L^\infty(0,T;\mbox{\boldmath $ V $}_0), 
\\
	\mbox{\boldmath $ v $} \in H^1(0,T;\mbox{\boldmath $ H $}_0)
	\quad \mbox{if } \tau >0, 
\\
	\mbox{\boldmath $ v $}^* {\pier{}:={}} (-\Delta v + {\pier\xi}, \partial_\nu v -\Delta_\Gamma v_\Gamma +{\pier\xi}_\Gamma  )
	\in L^2(0,T;\mbox{\boldmath $ H $}_0), 
\\
	\lambda \in L^2(0,T),  
\\
	 \mbox{\boldmath $ A$}_\tau \mbox{\boldmath $ v$}'(t)
	+ \mbox{\boldmath $ v$}^*(t)
	+ \lambda (t)\mbox{\boldmath $ w$}
	= P \Bigl( \mbox{\boldmath $ f$}(t) - \mbox{\boldmath $ \Pi $}_0 \bigl (\mbox{\boldmath $ v$}(t) \bigr ) \Bigr) 
	\quad \mbox{in } \mbox{\boldmath $ H$}_0,
	\ \mbox{for a.a.\ } t \in (0,T),
\\
	\mbox{\boldmath $ v $}^*(t) \in 
	\partial \varphi \bigl (\mbox{\boldmath $ v$}(t) \bigr ) 
	\quad \mbox{in } \mbox{\boldmath $ H$}_0,
	\ \mbox{for a.a.\ } t \in (0,T),  
\\
	\lambda (t) \mbox{\boldmath $ w $} \in 
	\partial I_{\mbox{\boldmath \scriptsize $\overline{ K }$}} 
	\bigl (\mbox{\boldmath $ v$}(t) \bigr ) 
	\quad \mbox{in } \mbox{\boldmath $ H$}_0, 
	\ \mbox{for a.a.\ } t \in (0,T), 
\\
	\mbox{\boldmath $ v$}(0) = \mbox{\boldmath $ v$}_0 
	\quad \mbox{in } \mbox{\boldmath $ H$}_0.  
\end{gather*}
Moreover, let us point out that 
$$
	\mbox{\boldmath $ v$}^*(t)+\lambda (t) \mbox{\boldmath $ w$} \in \partial 
	( \varphi +I_{\mbox{\boldmath \scriptsize $ K$}})\bigl (\mbox{\boldmath $ v$}(t) \bigr )
	\quad \mbox{in } \mbox{\boldmath $ H$}_0, 
	\ \mbox{for a.a.\ } t \in (0,T).
$$
Therefore, {\pier it is clear} that $\mbox{\boldmath $ v$}$ is the solution of the {C}auchy problem expressed by 
\eqref{(E)}--\eqref{(I.C.)}. We note that although the solution $\mbox{\boldmath $ v$}$ 
of this problem is uniquely determined, the 
auxiliary quantities $\mbox{\boldmath $ v$}^*$ and $\lambda $ are 
not unique in general (cf.\ \cite[Remark~3.3]{CF14}, \cite[Remark~2]{FK13}). 

\section*{Acknowledgments}
\noindent

The authors wish to express their heartfelt gratitude to 
{\pier {\takeshi professors} {G}oro {A}kagi and {U}lisse {S}tefanelli}, 
who kindly gave them {\pier the opportunity of exchange visits}
supported by {\pier the JSPS--CNR bilateral joint research 
project} \emph{Innovative Variational Methods for Evolution Equations}.  
{\pier The present note also benefits from a partial support of the MIUR--PRIN 
Grant 2010A2TFX2 ``Calculus of variations'' and the GNAMPA (Gruppo 
Nazionale per l'Analisi Matematica, la Probabilit\`a e le loro Applicazioni) 
of INdAM (Istituto Nazionale di Alta Matematica) for PC.}

%Authors are deeply grateful to the anonymous referees for 
%reviewing the original manuscript and for many valuable comments that 
%helped to clarify and refine this paper. 


\begin{thebibliography}{99}

\bibitem{Aik93}T.\ {A}iki,
	\newblock Two-phase {S}tefan problems with dynamic boundary conditions, 
	\newblock Adv.\ Math.\ Sci.\ Appl., {\bf 2} (1993), 253--270.
	
\bibitem{Aik95}T.\ {A}iki,
	\newblock Multi-dimensional {S}tefan problems with dynamic boundary conditions, 
	\newblock Appl.\ Anal., {\bf 56} (1995), 71--94.
	
\bibitem{Aik96}T.\ {A}iki,
	\newblock Periodic stability of solutions to some degenerate parabolic equations with dynamic boundary conditions, 
	\newblock J.\ Math.\ Soc.\ Japan, {\bf 48} (1996), 37--59.
	
\bibitem{Bar10}V.\ {B}arbu, 
	\newblock {\it Nonlinear differential equations of 
	 monotone types in {B}anach spaces}, 
	\newblock Springer, London, 2010. 
	
\bibitem{Bre73}H.\ {B}r\'ezis, 
	\newblock {\it Op\'erateurs maximaux monotones et semi-groupes de contractions dans les especes de {H}ilbert}, 
	\newblock North-Holland, Amsterdam, 1973.
	
\bibitem{BG87}F.\ {B}rezzi and G.\ {G}ilardi, 
	\newblock Partial differential equations, 
	\newblock H.\ {K}ardestuncer and D.\ H.\ {N}orrie (Eds.), 
	{\it Finite element handbook}, 
	McGraw-Hill Book Co., New York, 1987, {\pier Part~1.}
	
\bibitem{CH58}J.\ W.\ {C}ahn and J.\ E.\ {H}illiard, 
	\newblock Free energy of a nonuniform system I. Interfacial free energy, 
	\newblock J.\ Chem.\ Phys., {\bf 2} (1958), 258--267.
	
\bibitem{CC13}L.\ {C}alatroni and P.\ {C}olli,
	\newblock Global solution to the {A}llen--{C}ahn equation with 
	 singular potentials and dynamic boundary conditions, 
	\newblock Nonlinear Anal., {\bf 79} (2013), 12--27.
	
\bibitem{CGM13}L.\ {C}herfils, S.\ {G}atti and A.\ {M}iranville, 
	\newblock A variational approach to a {C}ahn--{H}illiard model in a domain 
	with nonpermeable walls,
	\newblock J.\ Math.\ Sci.\ (N.Y.), {\bf 189} (2013), 604--636. 
	
\bibitem{CF14}P.\ {C}olli and T.\ {F}ukao, 
	\newblock {A}llen--{C}ahn equation with dynamic boundary 
	 conditions and mass constraints, 
	\newblock Preprint arXiv:1405.0116~[math.AP] (2014), 
	 pp.~1--23{\pier , to appear in Math.\ Methods Appl.\ Sci.}
		
\bibitem{CGS14}P.\ {C}olli, G.\ {G}ilardi and J.\ {S}prekels, 
	\newblock On the {C}ahn--{H}illiard equation with dynamic 
	boundary conditions and a dominating boundary potential, 
	\newblock {\pier J.\ Math.\ Anal.\ Appl.\ {\bf 419} (2014), 972--994.}
{\pier%
\bibitem{CGS14bis}P.\ {C}olli, G.\ {G}ilardi and J.\ {S}prekels, 
	\newblock A boundary control problem for the viscous {C}ahn--{H}illiard equation with dynamic boundary conditions, 
	\newblock Preprint arXiv:1407.3916~[math.AP] (2014), pp.~1-27.
}%
\bibitem{CV90}P.\ {C}olli and A.\ {V}isintin, 
	\newblock On a class of doubly nonlinear evolution equations, 
	\newblock Comm.\ Partial Differential Equations {\bf 15} (1990), 737--756.
	
\bibitem{EZ86}C.\ M.\ {E}lliott and S.\ {Z}heng, 
	\newblock On the {C}ahn--{H}illiard equation, Arch.\ Ration.\ Mech.\ Anal., 
	{\bf 96} (1986), 339--357.
	
\bibitem{FK13}T.\ {F}ukao and N.\ {K}enmochi,
	\newblock Abstract theory of variational inequalities and {L}agrange multipliers, 
	\newblock pp.\ 237--246 in {\it Discrete and continuous 
	dynamical systems, supplement 2013}, 2013. 
	
\bibitem{GMS09}G.\ {G}ilardi, A.\ {M}iranville and G.\ {S}chimperna, 
	\newblock On the {C}ahn--{H}illiard equation with irregular potentials and dynamic boundary conditions,
	\newblock Commun.\ Pure.\ Appl.\ Anal., {\bf 8} (2009), 881--912.
	
\bibitem{GMS10}G.\ {G}ilardi, A.\ {M}iranville and G.\ {S}chimperna, 
	\newblock Long-time behavior of the {C}ahn--{H}illiard equation with irregular potentials and dynamic boundary conditions,
	\newblock Chin.\ Ann.\ Math.\ Ser.\ B, {\bf 31} (2010), 679--712.
	
\bibitem{GM13}G.\ R.\ {G}oldstein and A.\ {M}iranville,
	\newblock A {C}ahn--{H}illiard--{G}urtin model with dynamic boundary conditions,
	\newblock Discrete Contin.\ Dyn.\ Syst.\ Ser.\ S, {\bf 6} (2013), 387--400.
	
\bibitem{Gri09}A.\ {G}rigor'yan,
	\newblock {\it Heat kernel and analysis on manifolds}, 
	\newblock American Mathematical Society, International Press, Boston, 2009.
	
\bibitem{Ken07}N.\ {K}enmochi,
	\newblock Monotonicity and compactness methods for nonlinear variational inequalities, 
	\newblock M.\ {C}hipot (Ed.), {\it Handbook of differential equations: 
	Stationary partial differential equations}, {\bf Vol.4}, 
	North-Holland, Amsterdam (2007), 203--298. 
	
\bibitem{KN96}N.\ {K}enmochi and M.\ {N}iezg\'odka, 
	\newblock Viscosity approach to modelling non-isothermal diffusive phase separation,
	\newblock Japan J.\ Indust.\ Appl.\ Math., {\bf 13} (1996), 135--169.
	
\bibitem{Kub12}M.\ {K}ubo,
	\newblock The {C}ahn--{H}illiard equation with time-dependent constraint, 
	\newblock Nonlinear Anal., {\bf 75} (2012), 5672--5685. 
{\pier%
\bibitem{LM}J.-L.\ {L}ions and E.\ {M}agenes,
	\newblock {\it Non-homogeneous boundary value problems and applications},
	\newblock Vol.~I, Springer, Berlin, 1972.
}%

\bibitem{RZ03}R.\ {R}acke and S.\ {Z}heng, 
	\newblock The {C}ahn--{H}illiard equation with dynamic boundary conditions, 
	\newblock Adv.\ Differential Equations, {\bf 8} (2003), 83--110.
	
\bibitem{Sim87}J.\ {S}imon, 
	\newblock Compact sets in the spaces $L^p(0,T;B)$, 
	\newblock Ann.\ {\pier Mat.\ Pura.\ Appl.~(4)}, {\bf 146} (1987), 65--96.
\end{thebibliography}
\end{document}